\documentclass{amsart}
\usepackage{amsthm}
\usepackage{pinlabel}
\usepackage{enumitem}
\usepackage{color}
\usepackage{algorithm}
\usepackage{subcaption}
\usepackage{graphicx}
\usepackage{multicol}
\newlist{steps}{enumerate}{1}
\setlist[steps, 1]{label = Step \arabic*:}
\newtheorem{theorem}{\bf Theorem}[section]

\theoremstyle{definition}

\newtheorem{example}[theorem]{\bf Example}

\newcommand{\rme}{\mathrm{e}}
\newcommand{\rmi}{\mathrm{i}}

\begin{document}

\title{Complex optical vortex knots}

\author{Benjamin \textsc{Bode}}
\address{Instituto de Ciencias Matemáticas (ICMAT), Consejo Superior de Investigaciones Científicas (CSIC), Campus Cantoblanco UAM, C/ Nicolás Cabrera, 13-15, 28049 Madrid, Spain}
\address{Departamento de Matemática Aplicada a la Ingeniería Industrial, ETSIDI, Universidad Politécnica de Madrid, Rda. de Valencia 3, 28012 Madrid, Spain}
\email{benjamin.bode@upm.es}





\begin{abstract}
The curves of zero intensity of a complex optical field can form knots and links: optical vortex knots. Both theoretical constructions and experiments have so far been restricted to the very small families of torus knots or lemniscate knots. Here we describe a mathematical construction that presumably allows us to generate optical vortices in the shape of any given knot or link. We support this claim by producing for every knot $K$ in the knot table up to 8 crossings a complex field $\Psi:\mathbb{R}^3\to\mathbb{C}$ that satisfies the paraxial wave equation and whose zeros have a connected component in the shape of $K$. These fields thus describe optical beams in the paraxial regime with knotted optical vortices that go far beyond previously known examples.
\end{abstract}

\maketitle

\section{Introduction}\label{sec:intro}
Since 1867, when Tait's experiments involving smoke rings inspired Lord Kelvin to speculate that atoms were vortex knots in the aether \cite{kelvin}, many connections between knot theory and physics have been found. Lord Kelvin's theory turned out to be wrong of course, but the idea that closed vortex lines in an ideal fluid never change their topology is still important. More generally, topological stability of geometric structures in physical systems is a desirable property, since it implies that its topological invariants are quantities that are conserved and do not change with time, opening the door to possible applications for communication or data storage. Knotted structures have since been studied in a wide range of areas of physics \cite{sutcliffe, colloid, kamien, ma, proment, irvine, bodeem, acoustic, berry}, both theoretically and experimentally.

In the real world, where fluids are not ideal, knots decay and can change their topology. The study of the sequences of knots, knot cascades, that can occur during such processes is an active area of research \cite{simone, cascades}, see also \cite{dna} for a study of knot cascades in DNA.

Static solutions that describe physical systems that contain a knot are stable in another sense: the topology does not change under sufficiently small perturbations. This makes topological structures interesting candidates for future devices of communication \cite{info}, information storage \cite{coding} or quantum computation \cite{tqc}.

This article is concerned with complex-valued optical fields, which satisfy the paraxial wave equation and which have knotted zeros. In analogy with fluids and superfluids, such knots are called optical vortex knots. They have been studied both theoretically \cite{bd1, bd2, marknature} and experimentally \cite{photon, charge, small, ultrasmall} over the last years. While the arguments from \cite{daniel} can be adapted to prove that every knot can (theoretically) arise as an optical vortex knot, only very few have been constructed, observed or simulated explicitly. Almost all studies focus on the small family of torus knots \cite{photon, small, ultrasmall, charge}, which are particularly symmetric, or in some cases on the family of lemniscate knots \cite{marknature, lemniscate}. In this paper we introduce a new mathematical construction of optical vortex knots, or rather a construction of the corresponding optical fields containing the knots, which presumably allows us to create any knot as an optical vortex knot. In particular, this method produces optical vortex knots for the first 36 knots in the knot table, all knots of up to 8 crossings. The construction can easily be performed on a standard personal computer.



The paraxial wave equation
\begin{equation}\label{eq:para}
\nabla_{\perp}^2\Psi+2\rmi k\frac{\partial \Psi}{\partial z}=0,
\end{equation}
where $\nabla_{\perp}^2=\partial_x^2+\partial_y^2$ and $k$ is the wave number (which we will set to 1), is an approximation of the Helmholtz equation (with $\mathbf{A}(x,y,z)=\Psi(x,y,z)\rme^{\rmi kz}$) that holds in regimes where
\begin{equation}
\left|\frac{\partial^2 \Psi}{\partial z^2}\right|<<\left|k\frac{\partial \Psi}{\partial z}\right|.
\end{equation}
This is another way of saying that the angle between the wave vector $\mathbf{k}$ and the optical axis is small. In particular, laser beams with propagation direction along the $z$-axis can be described by solutions to the paraxial wave equation, where the modulus of the complex field describes the intensity (brightness of the beam) at any given point, while its complex argument $\arg(\Psi)$ describes the phase at each point. Naturally, the argument is not well-defined at points where $\Psi=0$. The set of such points is generically 1-dimensional and in particular can form knotted and linked closed loops: \textit{optical vortex knots}.

Calling such a line a vortex knot is justified in that the argument of $\Psi$ completes a full $2\pi$-rotation along any small meridian around the knot, so that the gradient field $\nabla\arg(\Psi)$ circulates around the vortex knot in the same way that a fluid circulates around a vortex.

The solutions to Eq.~\eqref{eq:para} are static. Optical vortex knots have also been studied outside of the realm of the paraxial wave equation. We showed in earlier work \cite{quasi} that every link type can be constructed as a subset of a stable vortex link of an electromagnetic field, i.e., as the common zeros of an electromagnetic field that satisfies Maxwell's equations in vacuum for all time. In particular, as time evolves, the knotted nodal structure moves through space, but never changes its topology. These fields cannot be expected to be monochromatic, so that experiments with laser beams as in the case of solutions to Eq.~\eqref{eq:para} are not possible.

In 2010, Dennis et al. presented a method to produce solutions of the paraxial wave equation with knotted vortices \cite{marknature}. Their procedure can be summarised as follows. For every knot $K$ there exists a complex-valued polynomial $F:\mathbb{R}^3\to\mathbb{C}$ in three real variables $x$, $y$ and $z$ whose zeros form $K$. This follows for example from the Nash-Tognioli theorem \cite{RAG}. If we have an explicit expression for the polynomial, we can numerically propagate the function $F(x,y,0)$, i.e., find the solution $\Psi$ of the paraxial wave equation whose restriction to the ($z=0$)-plane is equal to $F(x,y,0)$. This propagation corresponds to a decomposition of the function $F(x,y,0)$ into Laguerre-Gaussian beams. A priori these polynomial beams do not describe physical fields, since $F(x,y,0)$ diverges as $R=\sqrt{x^2+y^2}$ goes to infinity. However, the polynomial solutions can be embedded in Gaussian beams if the beam width is sufficiently large \cite{marknature}, resulting in an optical vortex knot in an actual physical field.


This approach resulted in analytic solutions as well as experiments confirming optical vortices in the shape of torus knots and lemniscate knots \cite{marknature}. These families of knots are comparatively simple and very symmetric. At that point it seemed like the main difficulty in finding (theoretical) optical vortex knots was to find explicit expressions for the polynomial $F$. The reason why torus knots lend themselves to this approach is that it is known that the $(p,q)$-torus knot is the intersection of the set of zeros of a polynomial $f:\mathbb{C}^2\to\mathbb{C}$, $f(u,v)=u^p-v^q$, with the 3-sphere $S^3\subset\mathbb{C}^2$ \cite{brauner, milnor}, so that the desired polynomial can be obtained by composing an inverse stereographic projection with the function $f$ (which results in a rational function) and clearing the denominator. Finding the desired polynomial $F$ for a general knot or link is not a trivial problem.

In \cite{bodepoly} we presented an algorithmic construction that takes as its input a knot or link (in the form of a braid word) and that produces the desired polynomial $F$. However, propagating the $(z=0)$-slices of these polynomials did not always produce the desired results for complicated knots. In general, the zeros of the solution $\Psi$ of the paraxial wave equation did not form the knot that was given by the zeros of $F$. Even for the three-twist knot $5_2$, which is the simplest knot in the knot table that is neither a torus knot nor a lemniscate knot, the procedure failed. We have an analytic expression for the polynomial $F$, but have found no solution $\Psi$ that reproduces the topology. This meant that despite the algorithm from \cite{bodepoly} the family of knots that were accessible to theoretical or experimental study was limited to the families of torus and lemniscate knots. Furthermore, the reason why the propagation procedure produced the correct topology in some cases, but not in others remained a mystery. The same problem already occured for some torus knots in \cite{marknature}.

In joint work with Hirasawa \cite{mikami} the author developed a variation of the construction of polynomials. Again, we obtain for any given knot or link a corresponding polynomial map $F$. The main result of this article is that propagating the $(z=0)$-slices of these new polynomials seems to reproduce the desired vortex topologies. This allows us to find analytic solutions of the paraxial wave equation for any given knot in the knot table up to 8 crossings. 

In fact, we present two constructions of optical fields in this article. The first is the propagation to the constructed polynomial $F(x,y,0)$, exactly as in Dennis's work. The only difference is that we use the polynomial obtained from the algorithm in \cite{mikami}. As in \cite{marknature} these polynomial beams correspond to physical solutions once they are embedded inside Gaussian beams of sufficiently large beam width. The second construction propagates $F(x,y,0)\rme^{-(x^2+y^2)/(2w^2)}$, so that the resulting field automatically decays at infinity like a Gaussian beam of beam width $w$. We find that for sufficiently small $w>0$ the field has an optical vortex knot that is the mirror image of the knot obtained from the first construction (up to isotopy). Thus either of the two methods can be used to create any knot type, but for the same input knot they produce knots that differ by a mirror reflection.

We have no mathematical proof that these methods are guaranteed to work for every knot or link, but we provide some arguments in addition to the many new-found examples that indicate that this could be the case. 

Our construction can be understood in the wider context of holography. It establishes a close geometric relationship between behaviour of the field $\Psi$ on the $(z=0)$-plane and the topology of its zeros in 3-dimensional space. Such connections have attracted more and more interest over the last years, as the holographic principle in certain string theories has illustrated its potential importance far beyond the area of optics \cite{holo, susskind}.






The rest of the article is structured as follows. In Section~\ref{sec:poly} we review the construction of polynomials introduced in \cite{mikami} and compare it to the earlier method from \cite{bodepoly}. In Section~\ref{sec:prop} we explain in more detail the two propagation techniques that produce the optical fields. Throughout the article we illustrate the different steps of the construction using the example of the knot $7_2$, which previously has not been seen as an optical vortex knot. In Section~\ref{sec:table} we present optical vortex knots for any knot in the knot table of up to 8 crossings. We conclude with a remark on the possibility of applying our considerations about optical vortex knots to quantum vortex knots in Section~\ref{sec:quantum}. 

\textbf{Acknowledgements:} This research was supported by the European Union’s Horizon 2020 research and innovation programme through the Marie Sklodowska-Curie grant agreement 101023017. The author would like to thank Mark Dennis for helpful discussions and advice and Danica Sugic for advice on the computational aspects of the project.

\section{Two constructions of polynomials}\label{sec:poly}
\subsection{Braids}
Both the construction developed with Dennis \cite{bodepoly} and our new construction with Hirasawa \cite{mikami} are built on the idea of braids. We briefly present the most important facts on braids. A much more detailed account can be found in \cite{braids}.

A \textit{geometric braid} $B$ on $s$ strands is a collection of $s$ disjoint curves in $\mathbb{C}\times[0,2\pi]$, where we think of the complex planes as horizontal and the interval as vertical, parametrized by their height coordinate:
\begin{equation}\label{eq:braidpara}
\bigcup_{j=1}^s(z_j(t),t),\quad t\in[0,2\pi],
\end{equation}
where $z_j:[0,2\pi]\to\mathbb{C}$, $j=1,2,\ldots,s$, are smooth functions such that for every $j\in\{1,2,\ldots,s\}$ there is a $k\in\{1,2,\ldots,s\}$ with $z_j(0)=z_k(2\pi)$ and $z_j^{(m)}(0)=z_k^{(m)}(2\pi)$, where the latter denotes equalities between $m$th derivatives for all $m\in\mathbb{N}$. Since the curves are assumed to be disjoint, there are no intersections between different strands. Since each strand is parametrized by its height coordinate, no strand can loop back on itself. The condition for the $z_j$s at the endpoints of the interval guarantee that the $n$ intersections with the plane $\mathbb{C}\times\{0\}$ match the intersections with the plane $\mathbb{C}\times\{2\pi\}$ setwise.

Two geometric braids are called braid isotopic, if they are (smoothly) isotopic in $\mathbb{C}\times[0,2\pi]$, maintaining the braid property and fixing the endpoints of the strands throughout the isotopy. We refer to a braid isotopy class as a \textit{braid}.

The set of braid isotopy classes of geometric braids on $s$ strands forms the braid group $\mathbb{B}_s$ on $s$ strands, whose group operation is vertical stacking of braids (with a reparametrization of the vertical interval). The identity element can be represented by the trivial braid, whose strands are vertical lines, i.e. the functions $z_j(t)$, $j=1,2,\ldots,s$, are constants. 

Projecting the strands of a braid into the ($\text{Im}(z)=0$)-plane results in a set of curves in the plane that can cross each other. We can assume (up to braid isotopy) that there are only finitely many crossings, i.e. only finitely many values $t_*\in[0,2\pi]$, where two strands, parametrized by $z_j(t)$ and $z_k(t)$, satisfy $\text{Re}(z_j(t_*))=\text{Re}(z_k(t_*))$, and that all crossings are transverse. We choose the convention that at a crossing the strand with smaller $\text{Im}(z)$-coordinate is the overpassing strand. This is usually depicted by deleting a part of the undercrossing strand in a neighbourhood of the crossing. The projection of the braid into the plane along with the information of which strand is crossing over which at each crossing is called a braid diagram.

For all values of $t$ for which there is no crossing there is a well-defined order on the strands of the braid, defined via the value of $\text{Re}(z_j(t))$, $j=1,2,\ldots,s$. The first strand (for some fixed $t_*\in[0,2\pi]$) is the strand parametrized by $z_j(t)$ such that $\text{Re}(z_j(t_*))<\text{Re}(z_k(t_*))$ for all $k\neq j$. The second strand is the strand parametrized by $z_{j'}$ such that $\text{Re}(z_{j'}(t_*))<\text{Re}(z_k(t_*))$ for all $k\neq j,j'$. In the braid diagram this simply means that at each fixed height the strands are labeled with increasing numbers from left to right.

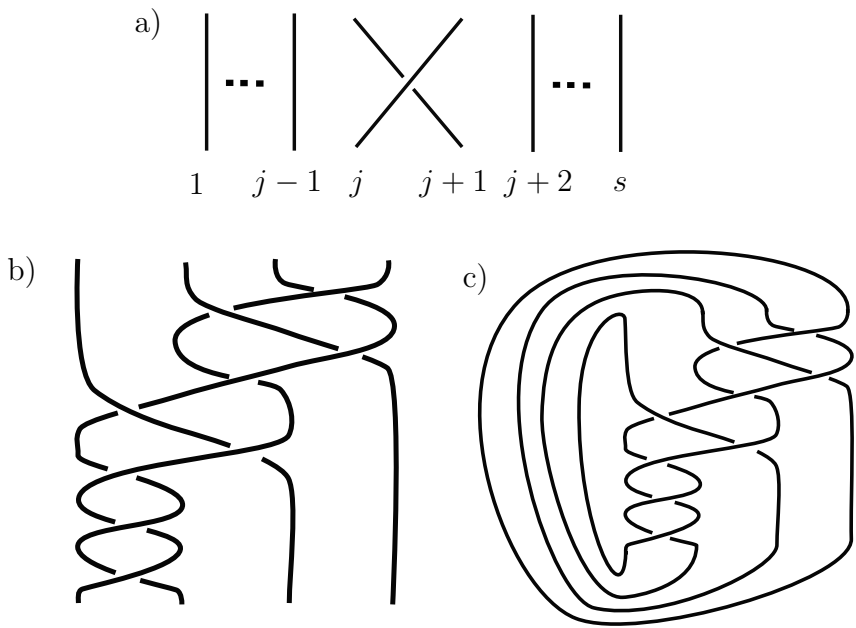
\begin{figure}[h]
\centering
\labellist
\Large
\pinlabel a) at -10 530
\pinlabel 1 at 150 0
\pinlabel $j-1$ at 450 0
\pinlabel $j$ at 680 0
\pinlabel $j+1$ at 1000 0
\pinlabel $j+2$ at 1280 0
\pinlabel $s$ at 1550 0 
\endlabellist
\includegraphics[height=2.5cm]{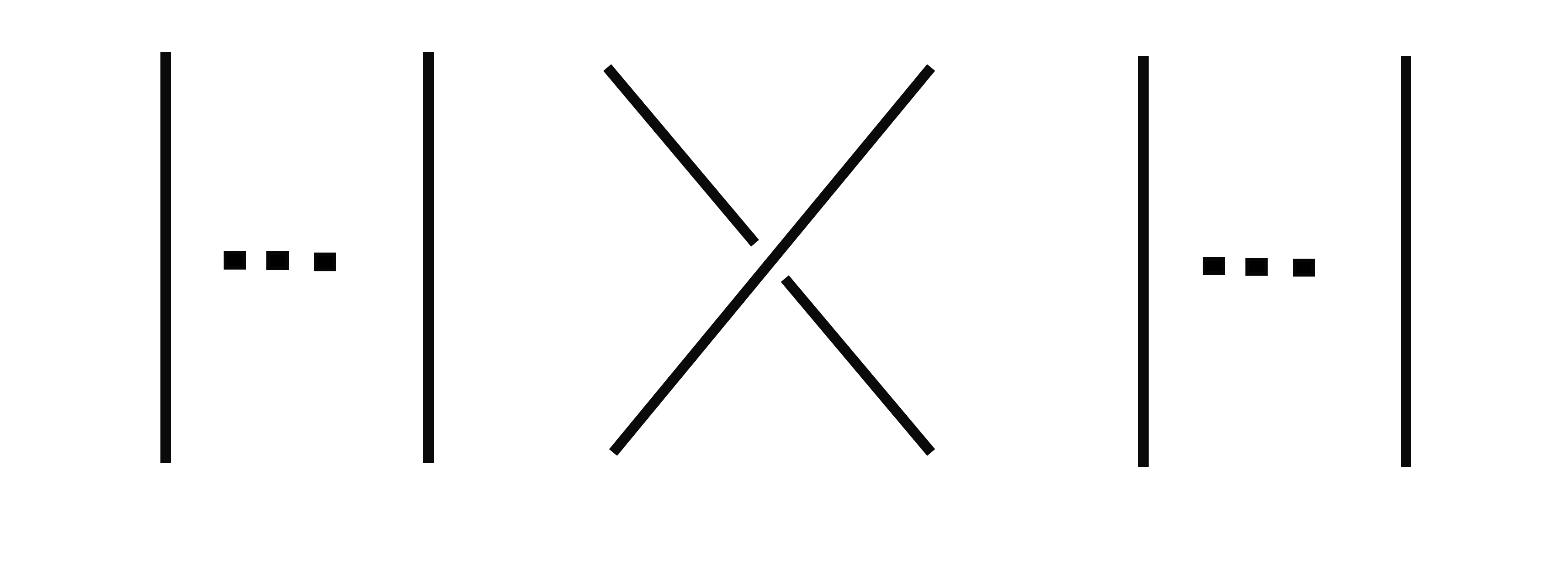}\\
\ \\
\ \\
\labellist
\Large
\pinlabel b) at 150 800
\endlabellist
\includegraphics[height=5cm]{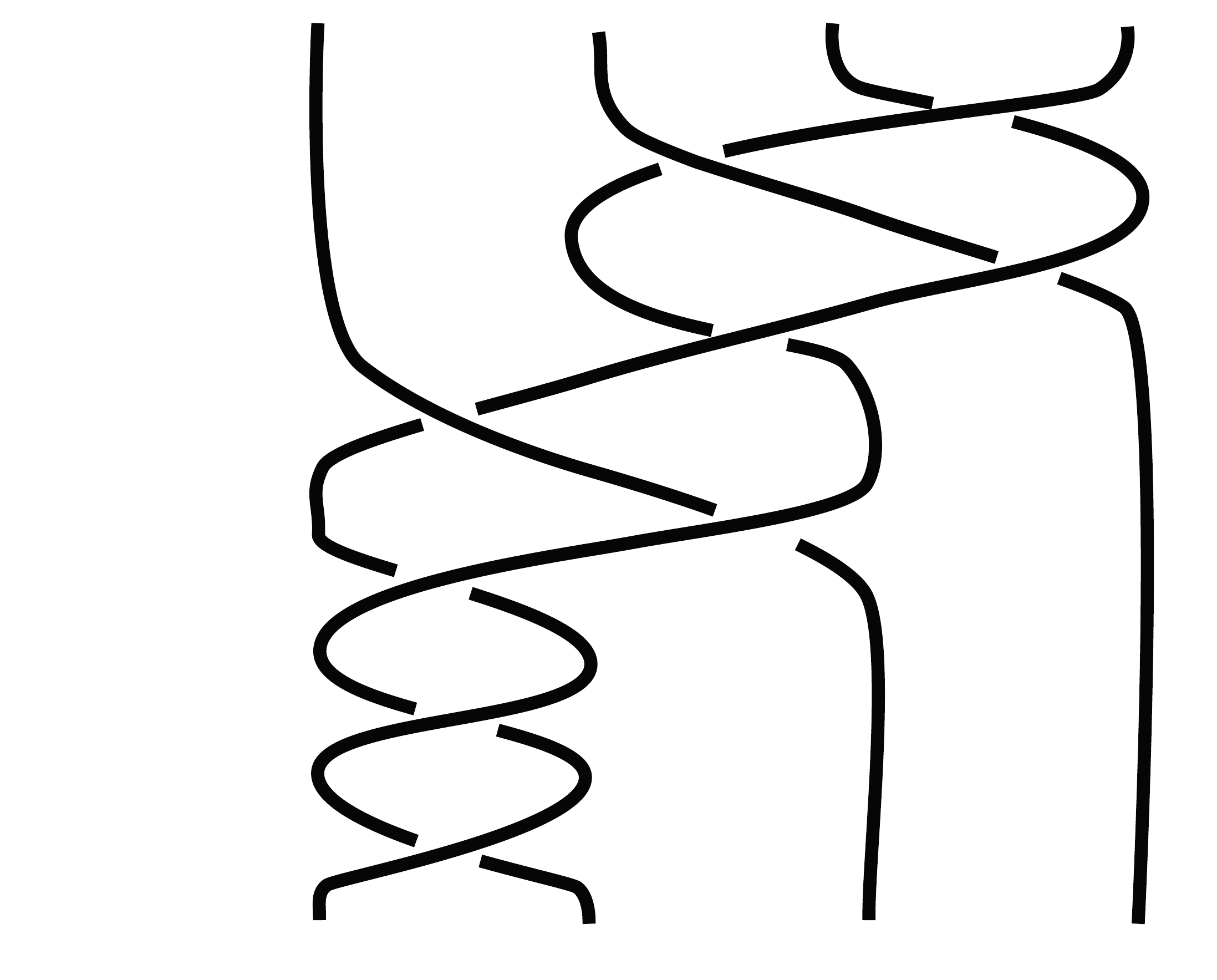}\qquad
\labellist
\Large
\pinlabel c) at 0 1100
\endlabellist
\includegraphics[height=5cm]{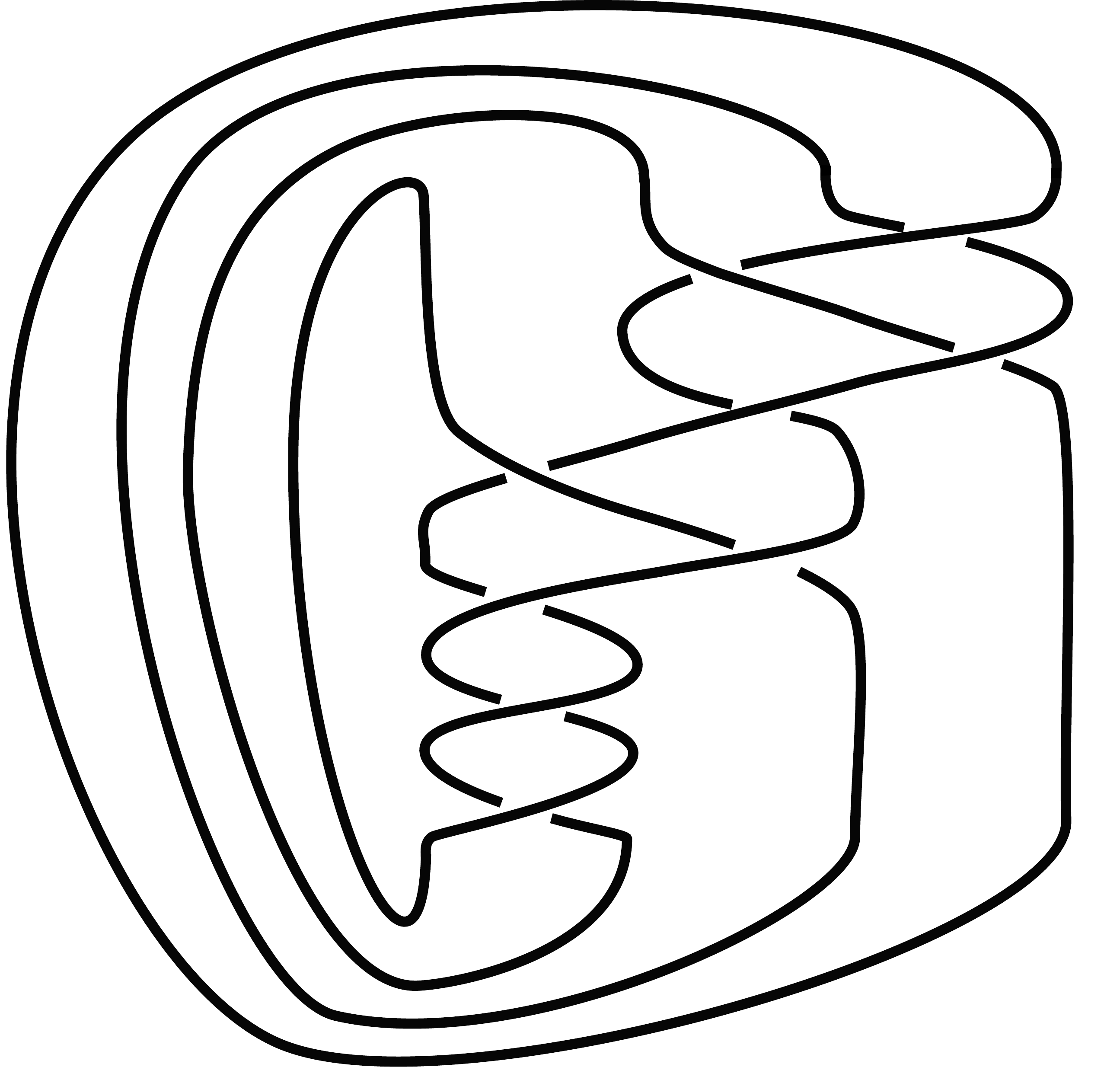}
\caption{a) The Artin generator $\sigma_j$. b) A braid diagram of the braid represented by the word $\sigma_1\sigma_1\sigma_1\sigma_2\sigma_1^{-1}\sigma_2\sigma_3\sigma_2^{-1}\sigma_3$. c) The closure of the braid, in this case the knot $7_2$. \label{fig:braids}}
\end{figure}

The most commonly used set of generators of $\mathbb{B}_s$ are the Artin generators $\sigma_j$, $j=1,2,\ldots,s-1$, which denote a positive half-twist between the $j$th strand and the $j+1$th strand. Figure~\ref{fig:braids}a) illustrates an Artin generator $\sigma_j$.

Since the endpoints of the strands at $t=0$ match the endpoints of the strands at $t=2\pi$ up to permuation, identifying the ($t=0$)-plane and the ($t=2\pi$)-plane results in a smooth link in the solid torus $\mathbb{C}\times S^1$. Using the standard (untwisted) embedding of the solid torus in the 3-sphere, we obtain a link in $S^3$, the closure of the braid. On the level of diagrams this closing procedure is realised by connecting the $s$ end points at the top of the braid to the $s$ end points at the bottom of the braid without introducing any new crossings. Figures~\ref{fig:braids}b) and c) show an example of a braid diagram and its closure, respectively.

The isotopy type of the braid closure in $S^3$ does not depend on the particular representative of a braid isotopy class, i.e., all braid isotopic geometric braids close to the same link type. Furthermore, every link is the closure of some braid \cite{alexander}.

\subsection{The construction by Bode-Dennis}\label{sec:mark}

Any geometric braid corresponds to a loop in the space of monic complex polynomials of degree $s$ and distinct roots. Another way of phrasing this is that the braid group on $s$ strands, the group that is formed by isotopy classes of braids on $s$ strands, is the fundamental group of the space of such polynomials. The identification sends a geometric braid as in Eq.~\eqref{eq:braidpara} to the loop of polynomials $g_t:\mathbb{C}\to\mathbb{C}$,
\begin{equation}
g_t(u):=\prod_{j=1}^s(u-z_j(t)).
\end{equation}

The height variable $t\in[0,2\pi]$ thus becomes the variable that parametrizes the loop. Note that the zeros of $g_t$ trace out exactly the given geometric braid as $t$ varies from $0$ to $2\pi$. We can interpret the loop of polynomials $g_t$ as one map $g:\mathbb{C}\times S^1\to\mathbb{C}$, $g(u,\rme^{\rmi t}):=g_t(u)$, whose zeros form the closed braid in the solid torus.

Note that for any $a>0$ the geometric braid
\begin{equation}\label{eq:lambda}
\bigcup_{j=1}^s(a z_j(t),t), t\in[0,2\pi],
\end{equation}
closes to the same link type as the original braid.

We think of $g$ as a holomorphic polynomial in the complex variable $u$, whose coefficients are $2\pi$-periodic functions of $t$. We can then approximate these coefficient functions by polynomials in $\rme^{\rmi t}$ and $\rme^{-\rmi t}$, without changing the braid type of the zeros. Substituting each instance of $\rme^{\rmi t}$ in $g$ by a second complex variable $v$ and every instance of $\rme^{-\rmi t}$ by its complex conjugates $\overline{v}$ results in a complex-valued polynomial $f$ in the complex variables $u$, $v$ and $\overline{v}$. Thus we can think of $f$ as a map from $\mathbb{C}^2$ to $\mathbb{C}$ (or from $\mathbb{R}^4$ to $\mathbb{R}^2$). By construction it is \textit{semiholomorphic}: holomorphic with respect to the variable $u$, but not necessarily with respect to $v$. Note also that $f(u,\rme^{\rmi t})=f(u,v)|_{v=\rme^{\rmi t}}=g(u,\rme^{\rmi t})$.

We showed in \cite{bodepoly} that for every geometric braid $B$ as above there is an $\varepsilon>0$ such that for all $a<\varepsilon$ the function $f$ that is built from Eq.~\eqref{eq:lambda} using the procedure above has the property that its zeros (or its \textit{vanishing set}) intersect $S^3\subset\mathbb{C}^2$ in the closure of $B$. Here $S^3$ explicitly denotes the 3-sphere of unit radius. This is not a statement about isolated singularities and 3-spheres of small radius a la Milnor \cite{milnor}.

Note that if $a$ is very small, then the complex coordinate of the braid parametrisation is very small. As a consequence the $u$-coordinate of the zeros of $g$ is very small. It follows that the intersection of the zeros $f$ and the 3-sphere lies in a tubular neighbourhood of $(0,\rme^{\rmi t})$, $t\in[0,2\pi]$.

Via composition with a stereographic projection map we obtain a complex-valued polynomial $F$ in three real variables $x$, $y$ and $z$, whose zeros form the closure of the given braid in $\mathbb{R}^3$. We obtain
\begin{equation}\label{eq:F}
F(x,y,z):=(x^2+y^2+z^2+1)^kf\left(\frac{x^2+y^2+z^2-1+2\rmi z}{x^2+y^2+z^2+1},\frac{2(x+\rmi y)}{x^2+y^2+z^2+1}\right),
\end{equation}
where $k$ is a sufficiently large natural number, so that the denominator on the right hand side is cancelled.

Since every link is the closure of a braid, we obtain a polynomial map for every link. In order to construct a given link $L$ in this way, we only have to find an explicit parametrisation of a geometric braid $B$ that closes to $L$. In \cite{bodepoly} this is done via a procedure involving trigonometric interpolation. Note that there is a lot of freedom in this construction. Every link is the closure of infinitely many braids and every braid admits infinitely many parametrizations, which all lead to different polynomials.

\subsection{The construction by Bode-Hirasawa}\label{sec:mikami}
We now describe the new construction of polynomials developed in joint work with Hirasawa, which will ultimately produce the correct vortex topology in solutions to the paraxial wave equation. First, consider a monic complex polynomial $p:\mathbb{C}\to\mathbb{C}$ of degree $s$. Assume that all of the roots $z_j$, $j=1,2,\ldots,s$, of $p$ are distinct and real and one of them is equal to zero. It then follows that all of the critical points $c_j$, $j=1,2,\ldots,s-1$, of $p$, i.e., the roots of $p'$, are also distinct and real. In fact, we may assume that the roots $z_j$ and the critical points $c_j$ are indexed according to their order in $\mathbb{R}$, so that $z_1<z_2<\ldots<z_s$ and $c_1<c_2<\ldots<c_{s-1}$. Then $c_j$ is the unique critical point in the interval between $z_j$ and $z_{j+1}$.

We write $v_j:=p(c_j)$, $j=1,2,\ldots,s-1$, for the critical values of $p$. After a small (real) deformation of $p$, we can also assume that the critical values are distinct and real. Note however, that their indexing does not contain any information about their order in $\mathbb{R}$.

By the Jordan Curve theorem every simple loop in $\mathbb{C}$ splits the complex plane into two connected components: a bounded component of its complement, which we will call the \textit{interior}, and an unbounded component, the exterior. Now take for every $j=1,2,\ldots,s-1$ a simple loop $\gamma_j(\chi)$, $\chi\in[0,2\pi]$, in $\mathbb{C}$ that is based at the origin, that is disjoint from all critical values and that contains exactly one critical value in its interior, namely $v_j$. Since it is a simple loop with $v_j$ in its interior, the loop $\gamma_j$ encircles $v_j$ exactly once. We require that $\gamma_j$ encircles $v_j$ in the counter-clockwise direction. In other words, the winding number of $\gamma_j$ relative to $v_j$ is positive. We write $\gamma_j^{-1}$ for the inverse loop, the same path traversed in the opposite direction, so that it encircles $v_j$ clockwise. 

Suppose now that $\prod_{k=1}^\ell\sigma_{j_k}^{\varepsilon_k}$, with $\varepsilon_k\in\{\pm 1\}$, is a braid word representing the braid $B$ with $s$ strands and $\ell$ crossings that we want to construct. Write $\Gamma(t)$, $t\in[0,2\pi]$, for the loop in $\mathbb{C}$ that is the concantenation of the loops $\gamma_{j_k}^{\varepsilon_k}$. In parametric form this reads
\begin{equation}\label{eq:loop}
\Gamma(t):=\gamma_{j_k}^{\varepsilon_k}(\ell t-2(k-1)\pi) \qquad \text{ if }t\in[2(k-1)\pi/\ell,2k\pi/\ell],\ k\in\{1,2,\ldots,\ell\}.
\end{equation}

We now define the loop of polynomials $g_t(u):=p(u)-\Gamma(t)$. Again, we can interpret this loop of polynomials as one map $g:\mathbb{C}\times S^1$, $g(u,\rme^{\rmi t})=g_t(u)$. It was shown in \cite{mikami} that the roots of $g$ form (up to braid isotopy) the desired (closed) braid $B$. As above, $g$ is a polynomial in the complex variable $u$ and its coefficients are $2\pi$-periodic functions of  $t$. In this case however, the only coefficient that depends on $t$ is the constant term $-\Gamma(t)$. Note that we can approximate $\Gamma(t)$ arbitrarily well by trigonometric polynomials, so that we can then apply the same procedure to the resulting function $g$ as in the previous subsection. We thus obtain a semiholomorphic polynomial $f$, which after an appropriate rescaling of the roots $z_j(t)$ by $a>0$ has the desired nodal set topology on $S^3$. We then obtain $F$ from $f$ as in the previous construction.


\section{Propagation}\label{sec:prop}

\subsection{Polynomial beams}
Let $K$ be any knot or link. Suppose that we have a polynomial $F:\mathbb{R}^3\to\mathbb{C}$ whose zeros form $K$. This polynomial could for example be the result of either of the two constructions from Section~\ref{sec:poly}.

Let $\Psi:\mathbb{R}^3\to\mathbb{C}$ denote the unique solution of the paraxial wave equation that agrees with $F$ on the ($z=0$)-plane, i.e., $\Psi(x,y,0)=F(x,y,0)$. We can find an analytic expression of $\Psi$ by propagating $F(x,y,0)$. This technique is also used in \cite{marknature}. A good summary of the background can be found in \cite{dani}.

In practice the propagation of $F(x,y,0)$ is achieved as follows. First we change coordinates to cylindrical coordinates $(R,\varphi,z)$. For tuples $(n,\ell)$ with $n-\ell$ even, the function
\begin{equation}
P_{n,\ell}(R,\varphi,z):=R^{|\ell|}\rme^{\rmi \ell\varphi}\tfrac{n-|\ell|}{2}!\left(2\rmi z\right)^{\tfrac{n-|\ell|}{2}}L_{\tfrac{n-|\ell|}{2},\ell}\left(\frac{R^2}{-2\rmi z}\right)
\end{equation}
is the unique solution of the paraxial wave equation that is equal to $R^n\rme^{\rmi \ell}$ in the $(z=0)$-plane. Here $L_{n,\ell}$ is the generalised Laguerre polynomial defined recursively via
\begin{align}
L_{0,|\ell|}(x)&=1,\\
L_{1,|\ell|}(x)&=1+\ell-x,\\
L_{n+1,|\ell|}(x)&=\frac{(2n+1+|\ell|-x)L_{n,|\ell|}(x)-(n+|\ell|)L_{n-1,|\ell|}}{n+1}.
\end{align}

Note that by construction each monomial of $F(R,\varphi,0)$ is of the form $c_{n,\ell}R^n\rme^{\rmi\ell\varphi}$ with $n-\ell$ even, i.e.,
\begin{equation}
F(R,\varphi,0)=\sum_{n,\ell}c_{n,\ell}R^n\rme^{\rmi \ell\varphi}
\end{equation}
with $c_{n,\ell}=0$ if $n-\ell$ is odd. Replacing every monomial $c_{n,\ell}R^n\rme^{\rmi\ell\varphi}$ by $c_{n,\ell}P_{n,\ell}(R,\varphi,z)$ is thus a complex-linear combination of solutions of the paraxial wave equation
\begin{equation}
\Psi(R,\varphi,z)=\sum_{n,\ell}c_{n,\ell}P_{n,\ell}(R,\varphi,z)
\end{equation} 
Since the partial differential equation is homogeneous and linear with constant coefficients, this implies that the obtained field $\Psi$ is a solution as well. Furthermore, it coincides with $F(R,\varphi,0)$ in the $(z=0)$-plane.

The solutions are polynomial beams and in particular they do not converge to 0 as $R$ goes to infinity. Therefore, they do not describe physical beams. However, as in \cite{marknature} they can be ``embedded'' in Gaussian beams with sufficiently large beam width, that is to say, if $\Psi$ has an optical vortex knot in the shape of $K$, then there exists a Gaussian beam, the result of propagating $\rme^{-R^2/(2w^2)}F(R,\varphi,0)$ with sufficiently large $w$, that also has an optical vortex knot in the shape of $K$. It is an important question how large $w$ has to be chosen in order to reproduce the desired knot type, i.e., what the minimal value of $w$ is that still leads to the correct topology. In particular, we would like to know if there exist values of $w$ that reproduce $K$ and that can be realised in an experiment.

Applying this procedure to the polynomial maps from Section~\ref{sec:mark} (or the polynomials for torus knots that go back to Milnor \cite{milnor} and Brauner \cite{brauner}) has produced several simple optical vortex knots \cite{marknature}, specifically torus knots and lemniscate knots, both of which are families of particularly symmetric knots. For example, using the polynomial $f(u,v)=u^2-v^3$, whose zeros form the trefoil knot, to define $F$ produces a complex field $\Psi$ whose zeros also form a trefoil knot. However, for polynomials that were constructed for other knots using the method from Section~\ref{sec:mark}, such as $5_2$, this method failed to reproduce the knot. The zeros of $\Psi$ were not the desired $5_2$ knot. It was already stated in \cite{marknature} that it is not clear what causes $\Psi$ to have zeros of the same (or different) topology as the zeros of $F$.

Instead of Laguerre polynomials it is also possible to propagate the polynomial by using Bessel beams \cite{bessel}.

We applied the procedure outlined above to polynomials $F$ that were constructed using the method in Section~\ref{sec:mikami}. Again we find that the topology of the zeros of $\Psi$ is not necessarily that of the zeros of $F$. However, we observe that it is usually the first strand of the braid where problems, such as interference with other components of the nodal set, occur. In $\mathbb{R}^3$ this first strand corresponds to parts of the knot that are comparatively close to the $z$-axis.

We resolved this issue by constructing polynomials $F$ not for the desired knot $K$, but for a split link of two components, one of which is $K$ and the other is the unknot, a planar circle. The braid that is used as input in the algorithmic construction of the polynomial consists thus of a vertical first strand next to the desired braid. If $B=\prod_{j=1}^\ell\sigma_{i_j}^{\varepsilon_j}$ is a braid that closes to $K$, then we use $\tilde{B}=\prod_{j=1}^\ell\sigma_{i_j+1}^{\varepsilon_j}$ as an input. In particular, $\tilde{B}$ has one more strand than $B$.

With this method we generated for every knot $K$ in the knot table of up to 8 crossings an optical field $\Psi$ with an optical vortex knot in the shape of $K$. Sometimes the extra unknot can also still be seen as a component of the nodal set of $\Psi$ as an unknotted circle that winds around the $z$-axis once.

\begin{example}\label{ex:72}
The knot $7_2$ is neither a torus knot nor a lemniscate knot. It has never been created or observed as an optical vortex knot before. It is the closure of the braid 
\begin{equation}
B=\sigma_1\sigma_1\sigma_1\sigma_2\sigma_1^{-1}\sigma_2\sigma_3\sigma_2^{-1}\sigma_3.
\end{equation} We thus need to apply the construction by Bode-Hirasawa to the braid 
\begin{equation}
\tilde{B}=\sigma_2\sigma_2\sigma_2\sigma_3\sigma_2^{-1}\sigma_3\sigma_4\sigma_3^{-1}\sigma_4,
\end{equation} 
whose closure consists of a copy of $7_2$ and an unknot, see Figure~\ref{fig:Btilde}.

\begin{figure}[h]
\centering
\labellist
\Large
\pinlabel a) at 150 800
\endlabellist
\includegraphics[height=4.5cm]{braid_diagram}\qquad
\labellist
\Large
\pinlabel b) at 50 800
\endlabellist
\includegraphics[height=4.5cm]{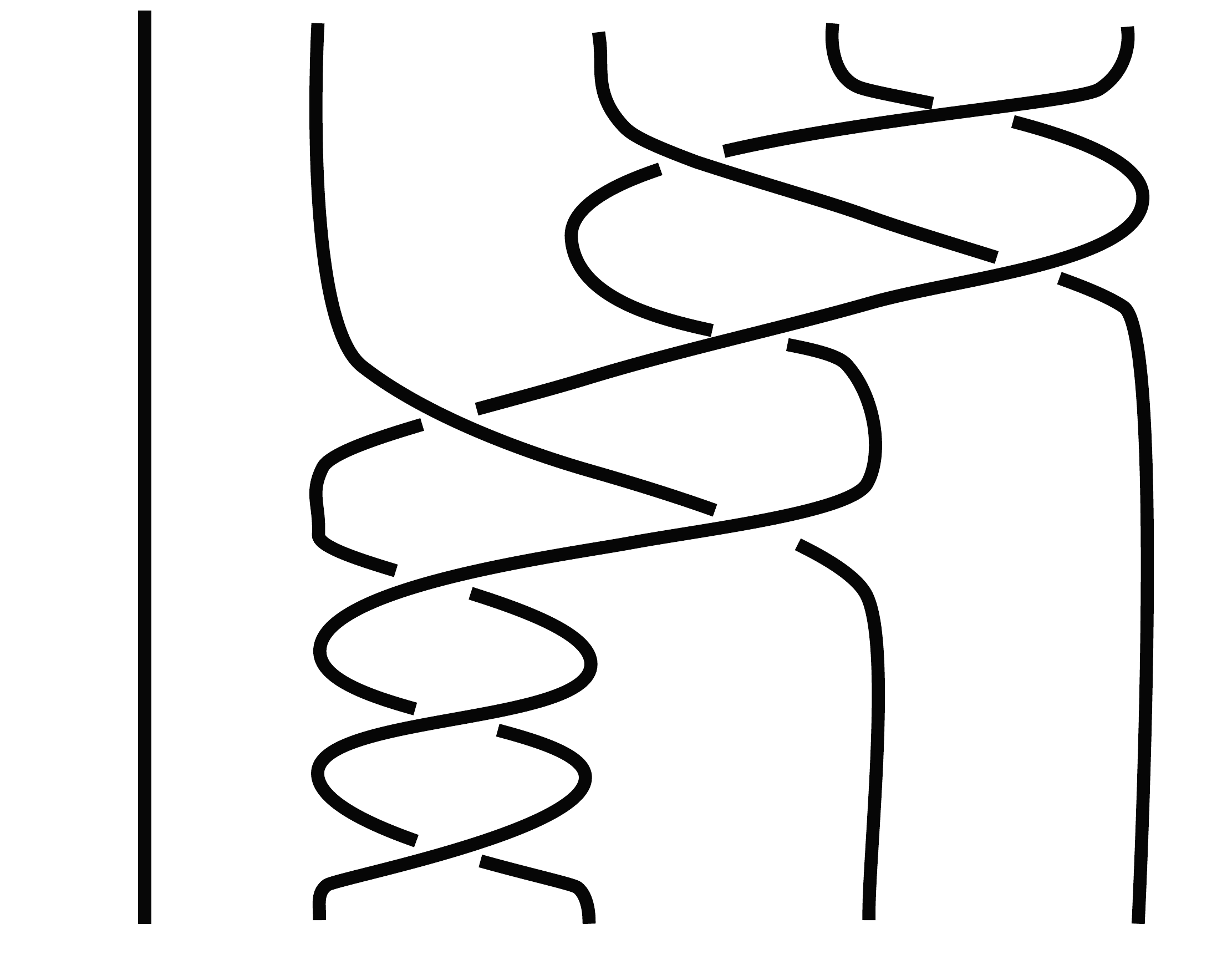}
\caption{a) A diagram of the braid $B$. b) A diagram of the braid $\tilde{B}$.\label{fig:Btilde}}
\end{figure}

Since it is a braid on 5 strands, we first need a polynomial in one variable of degree five, whose roots and critical points are real and simple. Take for example $p(u)=u (u - 1) (u + 1) (u - 2) (u + 2)$, whose roots are $z_1=-2,z_2=-1,z_3=0,z_4=1,z_5=2$ and whose critical points are $c_1=-1.6443,c_2=-0.543912,c_3=0.543912,c_4=1.6443$. Its critical values $v_j=p(v_j)$, $j=1,2,3,4$, are given by $v_1=3.63143,v_2=-1.4187,v_3=1.4187,v_4=3.63143$.

Next we need parametrisations of loops $\gamma_j$, $j=2,3,4$, in $\mathbb{C}$, that are based at the origin and wind around $v_j$ once. Since $\sigma_1$ does not appear in the braid word, we do not require a parametrisation of $\gamma_1$.
We use
\begin{align}
\gamma_2(\chi)&=1.3 (\cos(\chi) - 1) + 1.3 \rmi\sin(\chi), &\chi\in[0,2\pi],\nonumber\\
\gamma_3(\chi)&=1.3 (-\cos(\chi) + 1) - 1.3 \rmi\sin(\chi), &\chi\in[0,2\pi],\nonumber\\
\gamma_4(\chi)&=\begin{cases}2.2 (\cos(\chi)- 1) - 1.5\rmi \sin(2\chi) &\text{ if }\chi\in[0,\pi]\\
2.2 (\cos(\chi) - 1) +  1.5\rmi \sin(\chi) &\text{ if }\chi\in[\pi,2\pi].
\end{cases} &s\in[0,2\pi].
\end{align}

The loops are depicted in Figure~\ref{fig:basic72}. The same basic loops $\gamma_j$ can be used for any braid on 5 strands (that is, 4-stranded braids $B$ and the added vertical strand).

\begin{figure}[h]
\centering
\labellist
\Large
\pinlabel $v_1$ at 670 160
\pinlabel $v_2$ at 270 160
\pinlabel $v_3$ at 500 160
\pinlabel $v_4$ at 100 100
\pinlabel {{\color{blue}$\gamma_4$}} at 70 200
\pinlabel {{\color{red}$\gamma_2$}} at 250 250
\pinlabel {{\color{green}$\gamma_3$}} at 470 250
\endlabellist
\includegraphics[height=4cm]{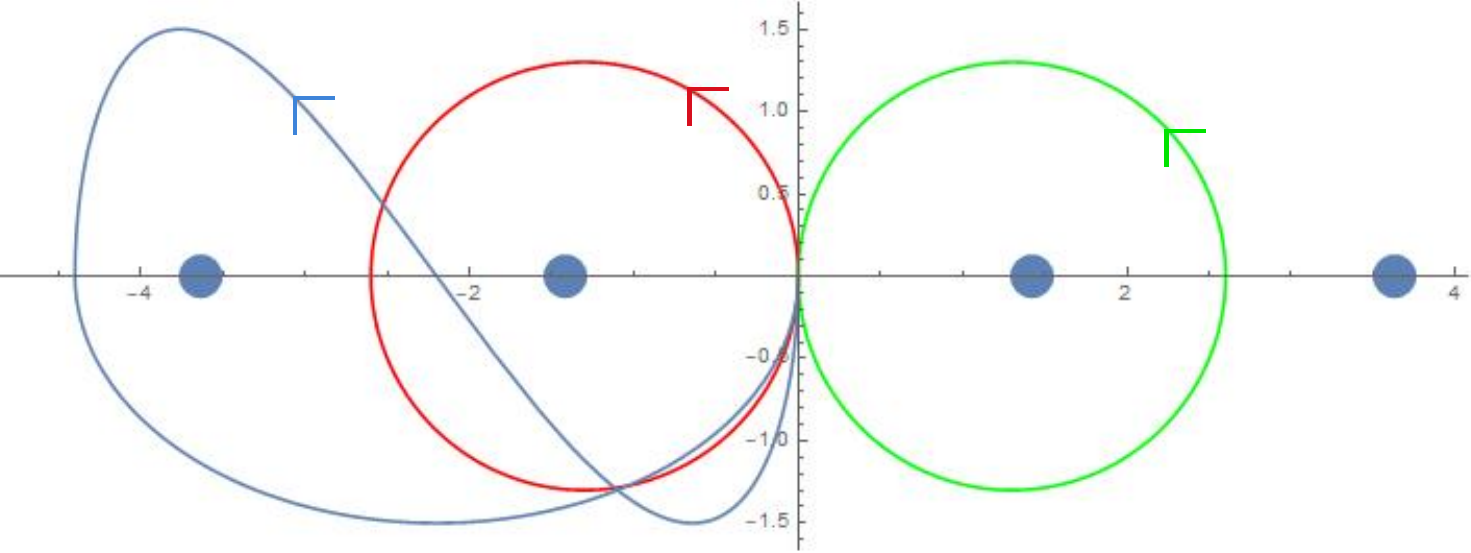}
\caption{The basic loops $\gamma_j(\chi)$, $j=2,3,4$, and the critical values $v_j$, $j=1,2,3,4$, of $p$.\label{fig:basic72}}
\end{figure}
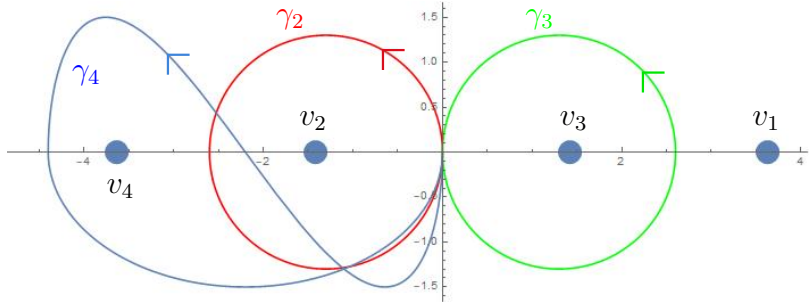

As in Eq.~\eqref{eq:loop} we obtain a loop $\Gamma(t)$ via an appropriate concatenation of these basics loops, which is determined by the braid word. 

We thus have a loop of polynomials $p(u)-\Gamma(t)$ whose zeros form the desired braid. We find that using the Fourier approximation $\Gamma_{trig}(t)$ of $\Gamma(t)$ of degree 20 instead of $\Gamma(t)$ in this polynomial expression produces a loop of polynomial $g_t(u):=p(u)-\Gamma_{trig}(t)$ with the same nodal topology. The roots of $g_t$ and the curve $(\Gamma_{trig}(t),t)$, $t\in[0,2\pi]$, in $\mathbb{C}\times[0,2\pi]$ are depicted in Figure~\ref{fig:critloop}.

\begin{figure}[h]
\centering
\labellist
\Large
\pinlabel a) at -30 400
\endlabellist
\includegraphics[height=6cm]{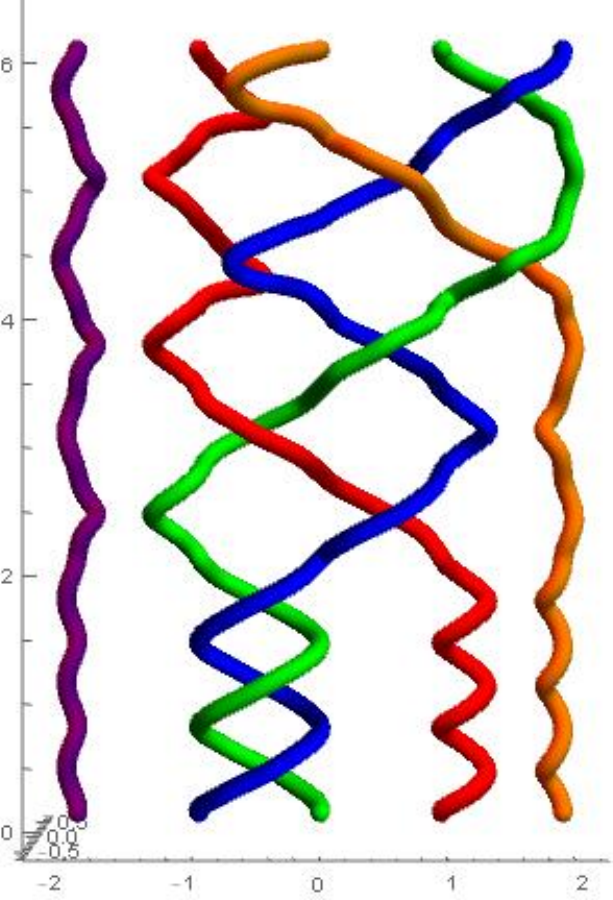}
\qquad
\labellist
\Large
\pinlabel b) at 20 500
\small
\pinlabel $v_4$ at 130 45
\pinlabel $v_2$ at 180 45
\pinlabel $v_3$ at 250 45
\pinlabel $v_1$ at 305 45
\endlabellist
\includegraphics[height=6cm]{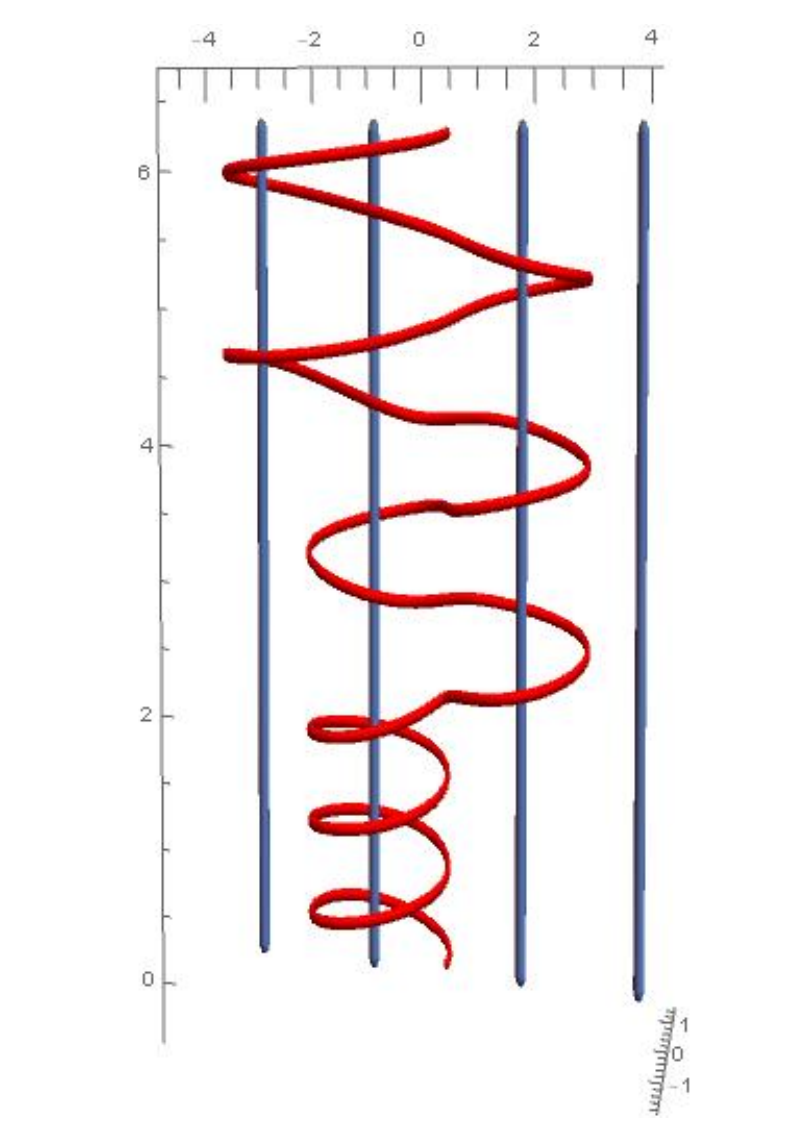}
\caption{a) The braid $\tilde{B}$ formed by the roots of $g_t$. b) The curve $(\Gamma_{trig}(t),t)$, $t\in[0,2\pi]$, in red and the curves $(v_i,t)$, $i=1,2,\ldots,4$, $t\in[0,2\pi]$ in blue, in $\mathbb{C}\times[0,2\pi]$. \label{fig:critloop}}
\end{figure}
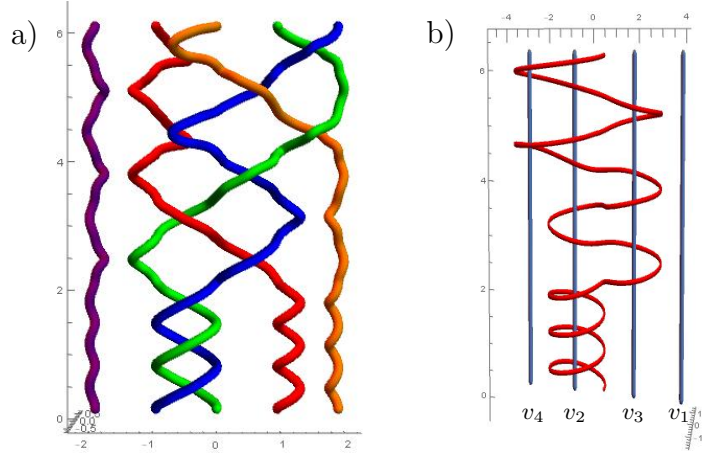
Figure~\ref{fig:critloop}a) also illustrates that the braid that is formed by the roots is not necessarily the desired braid on the nose. However, it only requires a straightforward isotopy, corresponding to some Reidemeister moves of type two \cite{reidemeister}, to see that it has the required topology.

We do not display all terms of $\Gamma_{trig}(t)$ here:
\begin{align}
\Gamma_{trig}(t)=&-0.624528-0.171786\rmi-(0.385917-0.0981406\rmi)\rme^{\rmi t}\nonumber\\
&+(0.0657862+0.226497\rmi)\rme^{2\rmi t}+\ldots-(0.0125504+0.0169264\rmi)\rme^{20\rmi t}\nonumber\\
&-(0.357661+0.134999\rmi)\rme^{-\rmi t}-(0.117603+0.199877\rmi)\rme^{-2\rmi t}\nonumber\\
&+\ldots+(0.0138012-0.00663464\rmi)\rme^{-20\rmi t}.
\end{align} 
A complete expression can be found in the corresponding mathematica file on the author's webpage \cite{webpage}.

Thus by definition
\begin{align}
g_t(u)=&u(u-a)(u+a)(u-2a)(u+2a)-\left[-0.624528-0.171786\rmi\right.\nonumber\\
&-(0.385917-0.0981406\rmi)\rme^{\rmi t}+(0.0657862+0.226497\rmi)\rme^{2\rmi t}\nonumber\\
&+\ldots-(0.0125504+0.0169264\rmi)\rme^{20\rmi t}\nonumber\\
&-(0.357661+0.134999\rmi)\rme^{-\rmi t}-(0.117603+0.199877\rmi)\rme^{-2\rmi t}\nonumber\\
&\left.+\ldots+(0.0138012-0.00663464\rmi)\rme^{-20\rmi t}\right].
\end{align}

and (with $a=\tfrac{1}{8}$)
\begin{align}
f(u,v)=&\tfrac{1}{8^5}8u\left(8u-1\right)\left(8u+1\right)\left(8u-2\right)\left(8u+2\right)-\left[-0.624528-0.171786\rmi\right.\nonumber\\
&-(0.385917-0.0981406\rmi)v+(0.0657862+0.226497\rmi)v^2\nonumber\\
&+\ldots-(0.0125504+0.0169264\rmi)v^{20}\nonumber\\
&-(0.357661+0.134999\rmi)\bar{v}-(0.117603+0.199877\rmi)\bar{v}^2\nonumber\\
&\left.+\ldots+(0.0138012-0.00663464\rmi)\bar{v}^{20}\right].
\end{align}

We define
\begin{equation}
F(x,y,z):=(x^2+y^2+z^2+1)^{20}f\left(\frac{x^2+y^2+z^2-1+2\rmi z}{x^2+y^2+z^2+1},\frac{2(x+\rmi y)}{x^2+y^2+z^2+1}\right),
\end{equation}
  
which in cylindrical coordinates reads
\begin{align}
F(R,\varphi,z)=\tfrac{1}{8^5}&\left(-30239.4-0.171786\rmi+(0.715322+0.269998\rmi)\rme^{-\rmi\varphi}R\right.\nonumber\\
&+(0.771834-0.196281\rmi)\rme^{\rmi\varphi}R-(292404-3.43572\rmi)R^2\nonumber\\
&+\ldots+(30240.6+0.171786\rmi)R^{40}+312384\rmi z+3.3753\times 10^{6}R^2 z+\ldots\nonumber\\
&\left.+(604812+3.43572\rmi)R^2z^{38}+312384\rmi z^{39}+(30240+0.17186\rmi)z^{40}\right).
\end{align}

Its restriction to the $(z=0)$-plane is
\begin{align}
F(R,\varphi,0)=\tfrac{1}{8^5}&\left(-30239.4-0.171786\rmi+(0.715322+0.269998\rmi)\rme^{-\rmi\varphi}R\right.\nonumber\\
&+(0.771834-0.196281\rmi)\rme^{\rmi\varphi}R-(292404-3.43572\rmi)R^2\nonumber\\
&\left.+\ldots+(30240.6+0.171786\rmi)R^{40}\right).
\end{align}

Propagating $F(R,\varphi,0)=\sum_{n,\ell}c_{n,\ell}R^n\rme^{\rmi\ell\varphi}$ yields the desired optical field
\begin{align}
\Psi(R,\varphi,0)=&\sum_{n,\ell}c_{n,\ell}P_{n,\ell}(R,\varphi,z)\nonumber\\
=&\tfrac{1}{8^5}\left(-30239.4-0.171786\rmi+(0.715322+0.269998\rmi)P_{1,-1}(R,\varphi,z)\right.\nonumber\\
&+(0.771834-0.196281\rmi)P_{1,1}(R,\varphi,z)-(292404-3.43572\rmi)P_{2,0}(R,\varphi,z)\nonumber\\
&\left.+\ldots+(30240.6+0.171786\rmi)P_{40,0}(R,\varphi,z)\right).
\end{align}  

A component of its nodal set is depicted in Figure~\ref{fig:72_poly}. As we can see it is the knot $7_2$, the closure of the braid $B$, except that all of the crossing signs have been switched. Strictly speaking it is therefore the mirror image of the desired braid closure. Usually (including in the knot tables) we do not distiguish between a knot and its mirror image even if they are not isotopic. Both are called the knot $7_2$. The figure was obtained by tracing the roots of $\Psi$, that is, we numerically find the zeros of $\Psi(R,\varphi,z)$ on the 840 planes $\varphi=2\pi j/840$, $j=0,1,2,\ldots,839$, and interpolate between pairs of nearest points.
  
\begin{figure}[h]
\centering
\includegraphics[height=6cm]{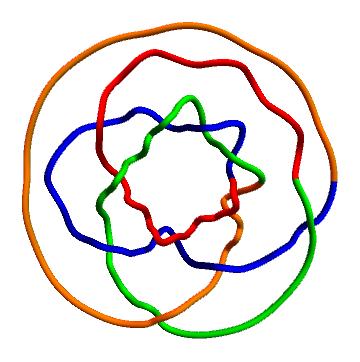}
\caption{A component of the nodal set of $\Psi$ is the optical vortex knot $7_2$.\label{fig:72_poly}}
\end{figure}
  
\end{example}

\subsection{Narrow Gaussian beams}

The polynomial beams from the previous subsection can be used to create optical vortex knots inside Gaussian beams as long as the beam width $w$ is sufficiently large. These beams are the result of propagating $F(R,\varphi,0)\rme^{-R^2/(2w^2)}$. Varying $w$ changes the topology of the zeros of the field. However, we find that for sufficiently small beam widths $w$ we again obtain an optical vortex knot. It is (up to isotopy) the mirror image of the knot that was obtained for large beam widths.

Instead of using the substitutions from the previous subsection, we propagate $F(R,\varphi,0)\rme^{-R^2/(2w^2)}$ as follows. Let 
\begin{equation}
F(R,\varphi,0)=\sum_{n,\ell}c_{n,\ell}R^n\rme^{\rmi \ell}.
\end{equation}
be the polynomial obtained from the method by Bode-Hirasawa, written in cylindrical coordinates. Note that (again) $c_{n,\ell}=0$ if $n-\ell$ is odd. The zeros of $F$ form the knot $K$. Note that
\begin{align}\label{eq:LGbeams}
Q_{n,\ell}(R,\varphi,z;w):=&\tfrac{n-|\ell|}{2}!\sum_{i=0}^{\tfrac{n-|\ell|}{2}}\binom{\tfrac{n-|\ell|}{2}+|\ell|}{\tfrac{n-|\ell|}{2} - i}R^{|\ell|}\rme^{\rmi\ell\varphi}\nonumber\\
&\times\frac{(1 - \rmi z/w^2)^
   n}{(1 + \rmi z/w^2)^{n + |\ell| + 1}} \rme^{\tfrac{-R^2}{2 w^2 (1 + \rmi z/w^2)}}w^{-|\ell|}L_{n,|\ell|}\left(\tfrac{R^2}{w^2 (1 + z^2/w^4)}\right)
\end{align}
is the unique solution of the paraxial wave equation that restricts to $\rme^{-R^2/(2w^2)}\left(\tfrac{R}{w}\right)^n\rme^{\rmi\ell\varphi}$ in the $(z=0)$-plane. As in the previous subsection $L_{n,|\ell|}$ is the generalised Laguerre polynomial. Thus 
\begin{equation}\label{eq:QNL}
\sum_{n,\ell}c_{n,\ell}Q_{n,\ell}(R,\varphi,z;1)
\end{equation}
is the unique solution that agrees with $\rme^{-R^2/2}F(R,\varphi,z)$ on the $(z=0)$-plane. This field does in general not reproduce the desired vortex knot.

Note that there is a scaling symmetry of the Laguerre-Gaussian beams in Eq.~\eqref{eq:LGbeams}. Since $Q_{n,\ell}(R,\varphi,z;w)=Q_{n,\ell}\left(\tfrac{R}{w},\varphi,\tfrac{z}{w^2};1\right)$, the propagation of $\rme^{-R^2/(2w^2)}F(R,\varphi,0)$ with beam width parameter $w=\mu$ is equal to the propagation of $\rme^{-R^2/2}F(\mu R,\varphi,0)$ with beam width parameter $w=1$. Clearly the scaling by a factor of $\mu$ or $\mu^2$ in the $R$- and $z$-coordinates does not affect the topology of the nodal set.



In other words, instead of keeping the polynomial $F$ fixed and varying the beam width $w$ (in order to find a sufficiently small value that reproduces the knot), we may keep the beam width $w=1$ fixed and vary the polynomial $F$ via a radial scaling. Instead of a knot that lies in a neighbourhood of the circle $R=1$, $z=0$, and a field with a very small beam width, we obtain a field with beam width 1 with a knot that lies in a neighbourhood of a planar circle of very large radius. The constructed field has the  analytic expression
\begin{equation}
\Psi_\mu(R,\varphi,z)=\sum_{n,\ell}c_{n,\ell}\mu^nQ_{n,\ell}\left(R,\varphi,z;1\right).
\end{equation}

Somewhat surprisingly at first, the zeros of the field $\Psi_\mu$ form the knot $m(K)$, i.e., the mirror image of the zeros of $F$, if $\mu$ is sufficiently small.

\begin{example}
We construct the knot $7_2$, (up to mirror image) the same knot as in Example~\ref{ex:72}. Recall that for this construction it is not necessary to add an unknotted component as in Example~\ref{ex:72}. We can apply the algorithm by Bode-Hirasawa directly to the braid
\begin{equation}
B=\sigma_1\sigma_1\sigma_1\sigma_2\sigma_1^{-1}\sigma_2\sigma_3\sigma_2^{-1}\sigma_3.
\end{equation}

Since the braid has 4 strands, we choose a polynomial of degree four as the starting point: $p(u)=u(u-1)(u+1)(u-1.8)$. Its critical points are $c_1=-0.621172,c_2=0.488653,c_3=1.48252$. Its critical values are $v_1=p(c_1)=-0.923653,v_2=p(c_2)=0.487784,v_3=p(c_3)=-0.5638$.

We have the following parametrisations for the basic loops $\gamma_j$, $j=1,2,3$:
\begin{align}
\gamma_1(\chi)&=0.53(\cos(\chi)-1)+\frac{\rmi}{2}(0.55\sin(\chi)-\sin(\chi)^2-0.2\cos(\chi)+0.2), &\chi\in[0,2\pi],\nonumber\\
\gamma_2(\chi)&=0.5(1-\rme^{\rmi \chi}), &\chi\in[0,2\pi],\nonumber\\
\gamma_3(\chi)&=0.375(\rme^{\rmi \chi}-1), &\chi\in[0,2\pi].
\end{align}
Finding the parametrizations of $\gamma_2$ and $\gamma_3$, which are the loops that wind around the critical values that are closest to 0, is straightforward, since the curves can be taken to be circles. Finding such a parametrisation for $\gamma_1$ requires a bit of trial and error. Note however, that once the parametrisations of these basic loops have been found, they can be used for any braid with the same number of strands (in this case 4). (A piecewise smooth parametrisation as in Example~\ref{ex:72} may also be used.)

The basic loops are displayed in Figure~\ref{fig:basic2}.

\begin{figure}[h]
\centering
\labellist
\Large
\pinlabel $v_3$ at 130 170
\pinlabel $v_1$ at 65 130
\pinlabel $v_2$ at 360 170
\pinlabel {{\color{blue}$\gamma_3$}} at 150 245
\pinlabel {{\color{red}$\gamma_1$}} at 20 220
\pinlabel {{\color{green}$\gamma_2$}} at 420 260
\endlabellist
\includegraphics[height=5cm]{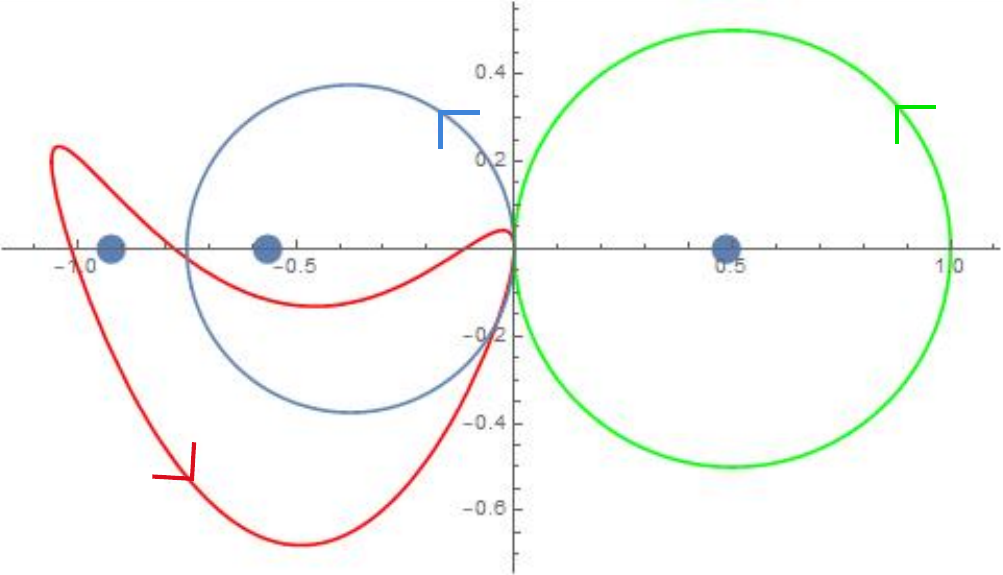}
\caption{The basic loops $\gamma_j(\chi)$, $j=1,2,3$, and the critical values $v_j$, $j=1,2,3$, of $p$.\label{fig:basic2}}
\end{figure}
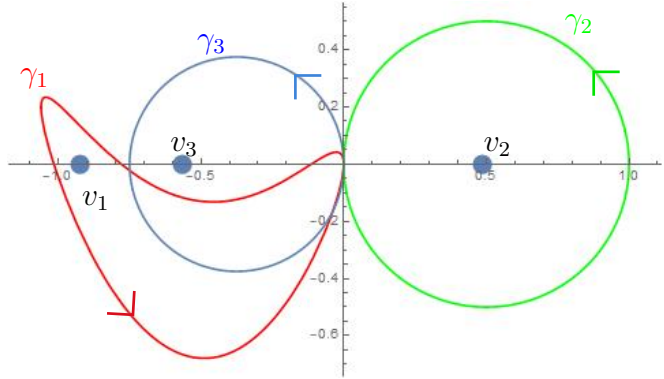

As in Example~\ref{ex:72} we combine the basic loops using Eq.~\eqref{eq:loop}, which results in a loop $\Gamma(t)$ such that the zeros of $p(u)-\Gamma(t)$ form the desired braid, a property that is shared by $g_t(u):=p(u)-\Gamma_{trig}(t)$, where $\Gamma_{trig}$ is the Fourier approximation of $\Gamma(t)$ of order 20.

The terms of lowest and highest order of $\Gamma_{trig}$ are
\begin{align}
\Gamma_{trig}(t)=&-0.150755-0.060854\rmi-(0.163382-0.112652\rmi)\rme^{\rmi t}\nonumber\\
&-(0.0147852-0.0906941\rmi)\rme^{2\rmi t}+\ldots+(0.00741904+0.00233376\rmi)\rme^{20\rmi t}\nonumber\\
&-(0.0448+0.129337\rmi)\rme^{-\rmi t}-(0.024861+0.11133\rmi)\rme^{-2\rmi t}\nonumber\\
&+\ldots-(0.00646206-0.00721021\rmi)\rme^{-20\rmi t}.
\end{align} 
For a complete expression for $\Gamma_{trig}$ we point the reader to the corresponding mathematica file on the author's webpage \cite{webpage}.

\begin{figure}[h]
\centering
\labellist
\Large
\pinlabel a) at -30 400
\endlabellist
\includegraphics[height=6cm]{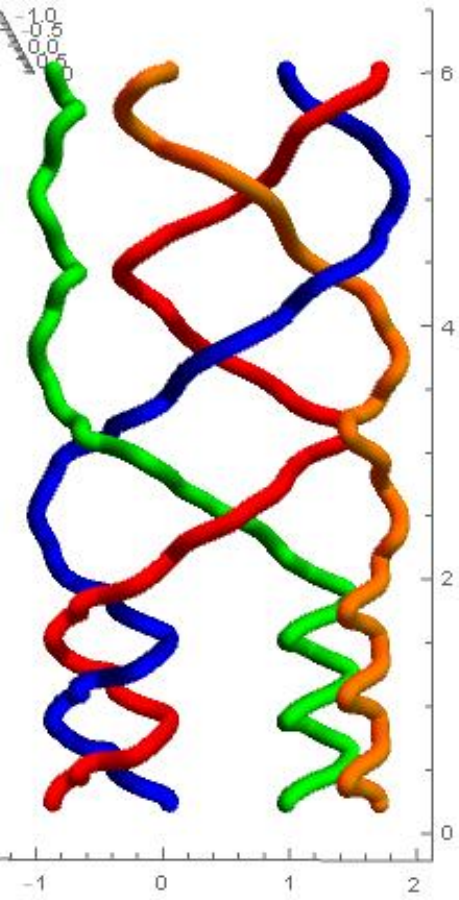}
\qquad
\labellist
\Large
\pinlabel b) at -30 400
\small
\pinlabel $v_1$ at 45 30
\pinlabel $v_3$ at 80 30
\pinlabel $v_2$ at 170 30
\endlabellist
\includegraphics[height=6cm]{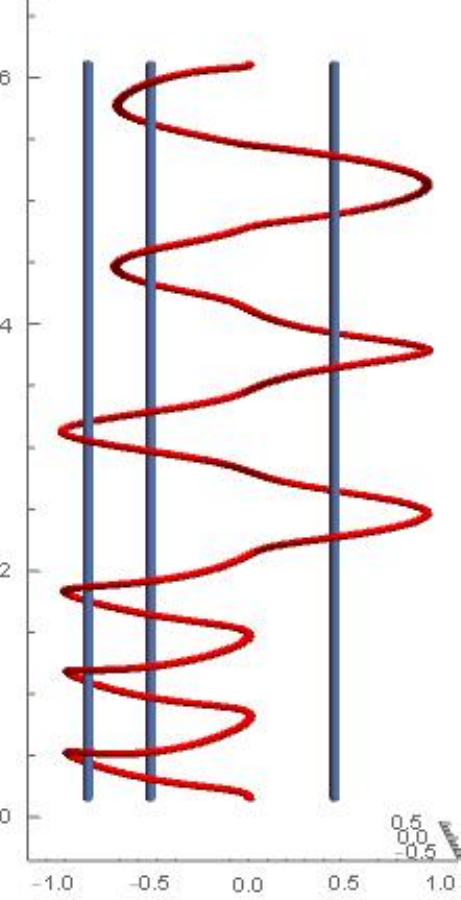}
\caption{a) The roots of the polynomial $g_t$ trace out the braid $B$. b) The curve $(\Gamma_{trig}(t),t)$, $t\in[0,2\pi]$, in red and the curves $(v_i,t)$, $i=1,2,\ldots,3$, $t\in[0,2\pi]$ in blue, in $\mathbb{C}\times[0,2\pi]$.\label{fig:72braid_roots}}
\end{figure}
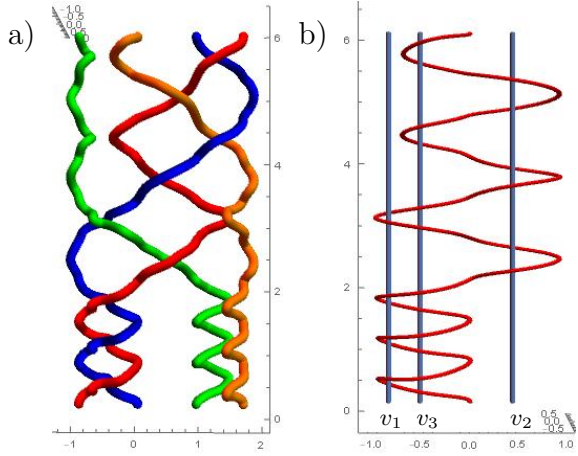

Thus by definition
\begin{align}
g_t(u)=&u(u-a)(u+a)(u-1.8a)-\left[-0.150755-0.060854\rmi\right.\nonumber\\
&-(0.163382-0.112652\rmi)\rme^{\rmi t}-(0.0147852-0.0906941\rmi)\rme^{2\rmi t}\nonumber\\
&+\ldots+(0.00741904+0.00233376\rmi)\rme^{20\rmi t}\nonumber\\
&-(0.0448+0.129337\rmi)\rme^{-\rmi t}-(0.024861+0.11133\rmi)\rme^{-2\rmi t}\nonumber\\
&\left.+\ldots-(0.00646206-0.00721021\rmi)\rme^{-20\rmi t}\right]
\end{align}

and (with $a=\tfrac{1}{12}$)
\begin{align}
f(u,v)=&\tfrac{1}{12^4}12u\left(12u-1\right)\left(12u+1\right)\left(12u-1.8\right)-\left[-0.150755-0.060854\rmi\right.\nonumber\\
&-(0.163382-0.112652\rmi)v-(0.0147852-0.0906941\rmi)v^2\nonumber\\
&+\ldots+(0.00741904+0.00233376\rmi)v^{20}\nonumber\\
&-(0.0448+0.129337\rmi)\bar{v}-(0.024861+0.11133\rmi)\bar{v}^2\nonumber\\
&\left.+\ldots-(0.00646206-0.00721021\rmi)\bar{v}^{20}\right].
\end{align}

We define
\begin{equation}
F(x,y,z):=(x^2+y^2+z^2+1)^{20}f\left(\frac{x^2+y^2+z^2-1+2\rmi z}{x^2+y^2+z^2+1},\frac{2(x+\rmi y)}{x^2+y^2+z^2+1}\right),
\end{equation}
  
which in cylindrical coordinates reads
\begin{align}
F(R,\varphi,z)=\tfrac{1}{12^4}&\left(23681+0.060854\rmi+(0.0896+0.258674\rmi)\rme^{-\rmi\varphi}R\right.\nonumber\\
&+(0.326764-0.225305\rmi)\rme^{\rmi\varphi}R+(289688+1.21708\rmi)R^2+\ldots\nonumber\\
&+(17503.4+0.060854\rmi)R^{40}-183931\rmi z-2.42587\times10^6\rmi R^2z\nonumber\\
&\left.+\ldots+(17503.4+0.060854\rmi)z^{40}\right).
\end{align}

Its restriction to the $(z=0)$-plane is
\begin{align}
F(R,\varphi,0)=\tfrac{1}{12^4}&\left(23681+0.060854\rmi+(0.0896+0.258674\rmi)\rme^{-\rmi\varphi}R\right.\nonumber\\
&+(0.326764-0.225305\rmi)\rme^{\rmi\varphi}R+(289688+1.21708\rmi)R^2+\ldots\nonumber\\
&\left.+(17503.4+0.060854\rmi)R^{40}\right).
\end{align}

Propagating $F(\mu R,\varphi,0)=\sum_{n,\ell}c_{n,\ell}\mu^nR^n\rme^{\rmi\ell\varphi}$ with $\mu=\tfrac{1}{12}$ and $w=1$ yields the desired optical field
\begin{align}
\Psi_{\mu}(R,\varphi,z)&=\sum_{n,\ell}c_{n,\ell}\mu^nQ_{n,\ell}(R,\varphi,z;1)\nonumber\\
\tfrac{1}{12^4}&\left(23681+0.060854\rmi+(0.0896+0.258674\rmi)\mu Q_{1,-1}(R,\varphi,z;1)\right.\nonumber\\
&+(0.326764-0.225305\rmi)\mu Q_{1,1}(R,\varphi,z;1)\nonumber\\
&+(289688+1.21708\rmi)\mu^2Q_{2,0}(R,\varphi,z;1)\nonumber\\
&+\ldots+(17503.4+0.060854\rmi)\mu^{40}Q_{40,0}(R,\varphi,z;1).
\end{align}  

Note that the construction gives an immediate decomposition of $\Psi$ into Laguerre Gaussian beams $Q_{n,\ell}$.

A component of its nodal set is depicted in Figure~\ref{fig:72}. As we can see it is the knot $7_2$, the closure of the braid $B$. It is thus (up to isotopy) the mirror image of the vortex knot that we obtained from the polynomial beam in Example~\ref{ex:72}. The figure was again obtained by tracing the roots of $\Psi$, through 840 planes of constant azimuthal coordinate $\varphi=2\pi j/840$, $j=0,1,2,\ldots,839$.
  
\begin{figure}[h]
\centering
\includegraphics[height=6cm]{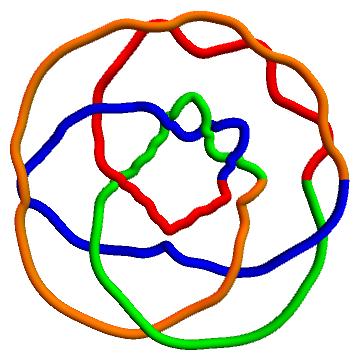}
\caption{A component of the nodal set of $\Psi$ is the optical vortex knot $7_2$.\label{fig:72}}
\end{figure}
\end{example}

This method seems a lot more reliable than the construction of polynomial beams in the sense that it is a lot easier to find values of $a$ and $\mu$ that lead to the correct knot. Furthermore, it is not necessary to add extra unknotted components. However, the knots obtained with this method might be harder to create and measure experimentally. It is a feature of this construction that the beam width is very small (compared to the size of the knot). Typically $w$ is chosen around $1/20$, while the knot lies close to circle of unit radius in the $(z=0)$-plane. Alternatively, as explained above, we can consider the beam as a solution with width $w=1$ and the knot lies close to the circle of radius $\mu^{-1}\approx 20$ in the $(z=0)$-plane. Therefore, in any case, the knot lies in a region of the field where the intensity is already very low.

\subsection{Knot holography}
We now present an argument why we expect our methods to reproduce the desired nodal topology for any given knot (or link).

The construction of polynomial beams is very similar to the earlier constructions. There is no clear mathematical argument that explains when this approach works and when it fails, but the general understanding is that there exist values for the parameter $a$ such that the optical field $\Psi$ is sufficiently close to the polynomial $F$ that they share the same nodal topology. One difficulty is finding such values. The values of $a$ that do not result in the correct knot type typically result in nodal configurations where a part of the knot has collided with other components of the zeros of $\Psi$. Shielding the knot from such a collision with an extra unknotted component seems to alleviate this problem.

The second method, the construction of Gaussian beams with small beam widths has a stronger mathematical foundation.
Note that, since the constant term $-\Gamma(t)$ is the only coefficient of $g$ that depends on $t$, we have $\tfrac{\partial g_t(u)}{\partial u}=p'(u)$ for all $t\in[0,2\pi]$. So the first derivative of $g$ with respect to the complex coordinate $u$ does not depend on $t$. In particular, the critical points of $g_t$, which are the roots of this derivative, are the same for every $t\in[0,2\pi]$, namely precisely $c_j$, $j=1,2,\ldots,s-1$, the critical points of $p$.

For any loop of monic polynomials $g_t$ of degree $s$ with distinct roots and distinct critical points, the roots form a braid on $s$ strands and the critical points form a braid on $s-1$ strands. In this case, the strands of the braid of critical points are parametrized by $\cup_{j=1}^{s-1}(c_j(t),t)$, but since here $c_j(t)=c_j$ for all $j$ and all $t\in[0,2\pi]$, the resulting braid of critical points is the trivial braid on $s-1$ strands. In \cite{mikami} the braid of critical points is called the \textit{saddle point braid} because it is formed by the set of saddle points of the surfaces of constant argument $\arg(g)$. Not only is the saddle point braid now the trivial braid, given by vertical strands, we also know that every $c_j$ is a real number. Therefore, the saddle point braid of $g_t$ lies in the plane of $\text{Im}(u)=0$.

Recall that under the stereographic projection map the set of points with $\text{Im}(u)=0$ is identified with the ($z=0$)-plane in $\mathbb{R}^3$, which is exactly the plane that is used to construct $\Psi$ by propagating. This is a first hint that important information on the function, namely, the values of the polynomial along the saddle point braid, can still be found in the plane that is used to construct the optical field.

To be precise, by construction we have $\frac{\partial f}{\partial u}(c_j,r\rme^{\rmi t})=0$ for all $j=1,2,\ldots,s-1$, $t\in[0,2\pi]$, $r> 0$. Let $V_1$ and $V_2$ be two linearly independent directions tangent to $S^3$ and orthogonal to $\partial_t$ using the round metric. Let $C$ denote the set of points on $S^3$ where the directional derivatives of $f$ in both directions $V_1$ and $V_2$ vanish. It follows from the same arguments as in \cite{bodepoly} that for smaller and smaller values of $a$ the approximation between $C$ and $(a c_j,\rme^{\rmi t})$ becomes better and better in the sense that if $(u_a(t),\rme^{\rmi t})$ is a local parametrisation of $C$ for a fixed parameter $a$, then $\tfrac{1}{a}u_a(t)-c_j$ converges to 0 as $a$ goes to $0$. In particular, $C$ consists of $s-1$ curves $C=\cup_{j=1}^{s-1}C_j$, such that $C_j$ is close to $(a c_j,\rme^{\rmi t})$. Thus for sufficiently small values of $a$ the curves $C$ lie very close to the ($\text{Im}(u)=0$)-plane.

Stereographic projection maps $\partial_t$ to the derivative $\partial_\varphi$ with respect to the angle $\varphi$ in cylindrical coordinates. Thus the image of $C$ under stereographic projection consists of a set of curves that are very close to the ($z=0$)-plane and that are zeros of two linearly independent directional derivatives of $F$, along directions that are orthogonal to $\partial_\varphi$.

The values of $F$ on $C$ are (up to an overall real factor) close approximations to the values that $g$ takes on the saddle point braid, which are exactly $v_j-\Gamma(t)$. We claim that these values essentially encode the topology of the zeros.

Note that 
\begin{align}
\frac{\partial \arg(g_t)}{\partial t}(c_j)&=\frac{\partial\text{Im}\log(g_t)}{\partial t}(c_j)\nonumber\\
&=\frac{\partial\text{Im}\log(\prod_{i=1}^s(u-z_i(t)))}{\partial t}(c_j)\nonumber\\
&=\sum_{i=1}^s\frac{\partial\text{Im}\log(u-z_i(t))}{\partial t}(c_j)\nonumber\\
&=\sum_{i=1}^s\text{Im}\left(\frac{-z_i'(t)}{c_j-z_i(t)}\right)\nonumber\\
&=\sum_{i=1}^s\frac{\text{Im}(-z_i'(t))\text{Re}(c_j-z_i(t))-\text{Re}(-z_i'(t))\text{Im}(c_j-z_i(t))}{|c_j-z_i(t)|^2}.
\end{align}

This last sum is thus exactly the total angular velocity of the $s$ roots $z_i(t)$ relative to $c_j$ as they move in the complex plane. Since each term is weighted by the inverse squared distance of the root $z_i(t)$ to $c_j$, the largest contribution to the sum comes from the roots that are closest to $c_j$.

Recall that $\Gamma(t)$ is a concatenation of the basic loops $\gamma_j(\chi)$, each of which has $0\in\mathbb{C}$ as its basepoint. Therefore $\Gamma\left(\tfrac{2\pi k}{\ell}\right)=0$ for all $k=0,1,2,\ldots,\ell$, where $\ell$ is the number of crossings in the desired braid. Then $g_{\tfrac{2\pi k}{\ell}}=p$ for all $k=0,1,2,\ldots,\ell$, which means that in regular intervals all roots return to their initial configuration along the real line (up to permutation). In particular, after each $2\pi/\ell$-interval all roots lie again in the $(\text{Im}(u)=0)$-plane.

Consider such an interval, for example $[0,2\pi/\ell]$. By construction $\Gamma(t)=\gamma_{j_1}^{\varepsilon_1}(\ell t)$ on this interval, where $j_1$ is the index of the first Artin generator in the desired braid word. Since $g_{\tfrac{2\pi}{\ell}}=g_{0}=p$ the image of the map $t\mapsto \arg(g_t(c_j))$ is a loop in $S^1$ for each $j=1,2,\ldots,s-1$. Each of these loops has a winding number, which is exactly the number of times that $\gamma_{j_1}^{\varepsilon_1}$ winds around $v_j=p(c_j)$. It is equal to $\varepsilon_1$ if $j=j_1$ and equal to zero if $j\neq j_1$. 

If the winding number is zero, this means that the net relative motion of the two roots closest to $c_j$ is zero, i.e., they do not interchange their position. Since this is a topological argument, we do not know the precise motion of the nearby roots. They could interchange their position an even number of times, but the vanishing of the winding number guarantees that the crossings generated in this way come in pairs with opposite signs and cancel each other.

If the winding number is 1, this means that the two roots closest to $c_j$ are both moving in a clockwise direction around $c_j$, each performing half of a full rotation around $c_j$, thereby generating a positive crossing $\sigma_j$. Similarly, a negative winding number implies counterclockwise motion and a negative crossing $\sigma_j^{-1}$. There are several geometric relations between the roots and the critical points of a complex polynomial that guarantee that the non-zero winding number has to correspond to a crossing. For example, it is not possible to have one root $z_j$ be stationary and the root $z_{j+1}$ wind once around the critical point $c_j$, since at some point $z_{j+1}$ would have to cross the real line between a critical point and a root, contradicting the principle that there needs to be a real critical point between any pair of real roots.

This is a somewhat vague interpretation of the argument why the roots of the constructed polynomials $g_t$ trace out the desired braid. The precise mathematical arguments are worked out elsewhere (for example, \cite{mikami} or \cite[Chapter 5]{phd}). The interpretation in terms of total angular velocity has not appeared before. The important insight is that $\arg(g_t(c_j))$, more specifically, their winding numbers on each $2\pi/\ell$-interval completely determine the topology of the roots. Since the critical points $c_j$ are real, they lie in the $(\text{Im}(u)=0)$-plane. The image of the corresponding curves $C_j$ under stereographic projection lies close to the $(z=0)$-plane, but not necessarily inside that plane. However, by continuity there are tubular neighbourhoods of the curves $C_j$ that intersect the $(z=0)$-plane in $s-1$ concentric annuli, such that along each $2\pi/\ell$-interval of the azimuthal coordinate $\varphi$ the phase of the field has the same winding behaviour as $\arg(g)$ along the saddle point braid.

\begin{figure}[h]
\centering
\includegraphics[height=5cm]{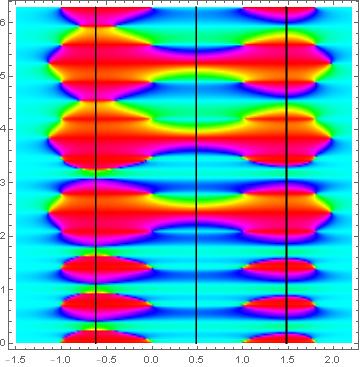}
\includegraphics[height=1cm]{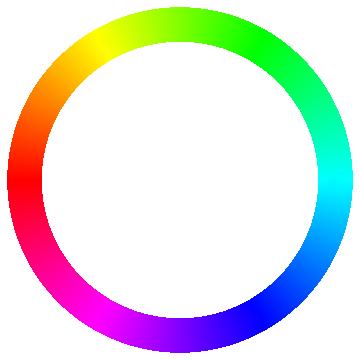}\qquad
\includegraphics[height=5cm]{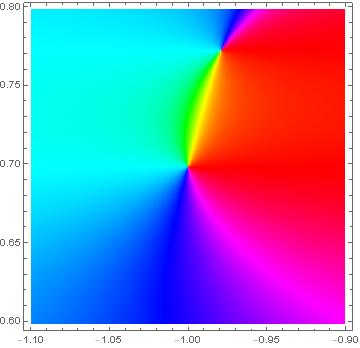}\\
\ \\
\labellist
\Large
\pinlabel a) at -240 750
\pinlabel b) at 170 750
\pinlabel c) at -10 410
\tiny 
\pinlabel $\text{Re}(u)$ at 70 455
\pinlabel $t$ at -220 760
\pinlabel $\text{Re}(u)$ at 500 455
\pinlabel $t$ at 190 760
\pinlabel 0 at 158 492
\pinlabel $\pi/2$ at 125 530
\pinlabel $\text{Re}(u)$ at 275 0
\pinlabel $t$ at 15 420
\endlabellist
\includegraphics[height=7cm]{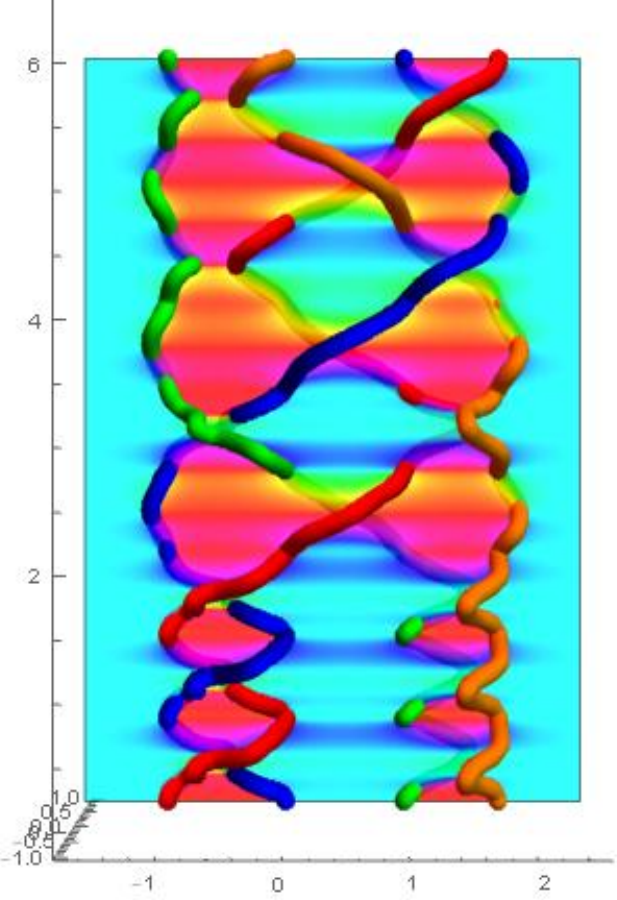}
\caption{a) The values of $\arg(g)$ on the $(\text{Im}(u)=0)$-plane. The vertical black lines are the saddle point braid, formed by $c_j\times[0,2\pi]$. b) Two intersection points between the braid $B$, formed by the roots of $g$, and the $(\text{Im}(u)=0)$-plane. Around each of them the phase rotates by $2\pi$. However, from further away it might appear like some colours (like green or yellow) do not appear. c) The braid $B$ and $\arg(g)$ on the $(\text{Im}(u)=0)$-plane. \label{fig:72imu0}}
\end{figure}
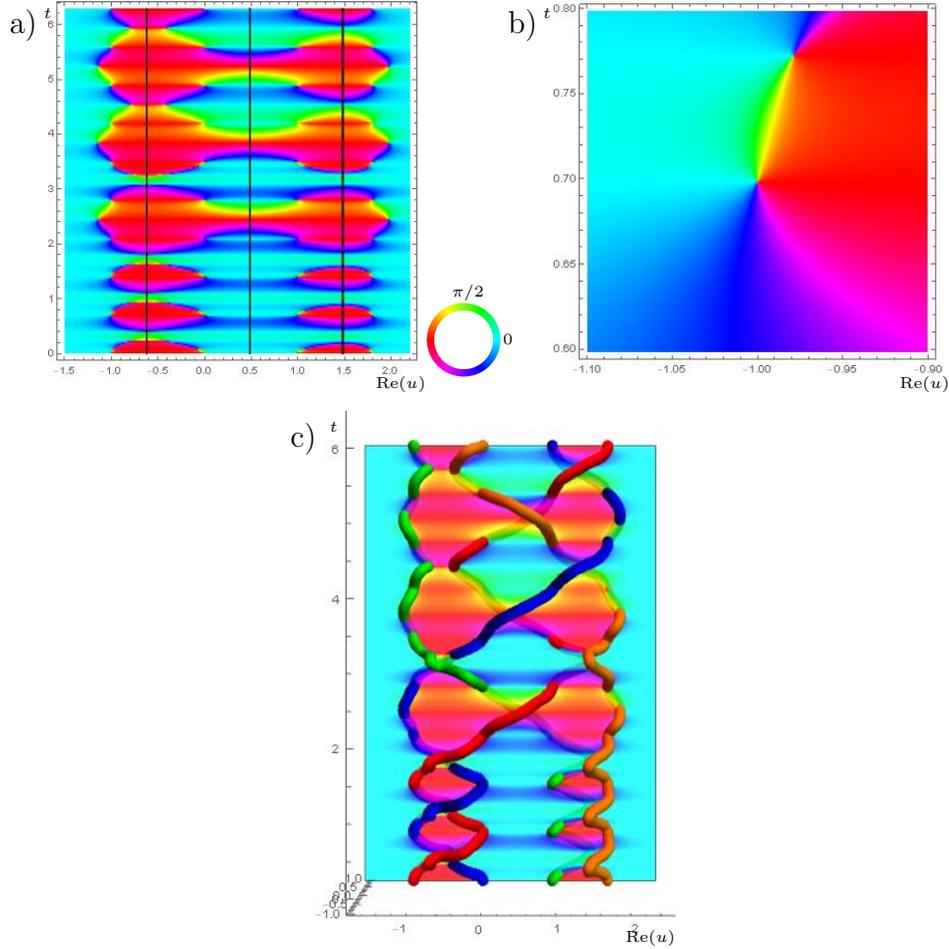

We illustrate our observations once again using the example of the knot $7_2$. Figure~\ref{fig:72imu0}a) shows a plot of the function $\arg(g_t)$ on the plane $\text{Im}(u)=0$. The black vertical lines are the curves $(c_j,t)$, $j=1,2,\ldots,s-1$, $t\in[0,2\pi]$ i.e., the saddle point braid, which is entirely contained in that plane. The points where the braid $B$, formed by the roots of $g_t$ intersects the plane can be identified, since around them the phase rotates by $2\pi$. In Figure~\ref{fig:72imu0}a) it might appear like there are points around which the phase rotates by some number that is not $2\pi$. In fact, these are simply two roots that are very close to each other, as shown in Figure~\ref{fig:72imu0}b). Figure~\ref{fig:72imu0}c) shows the braid that is formed by the roots of $g_t$ and the plot of $\arg(g_t)$ on the plane $\text{Im}(u)=0$. Note the regular intersections of the braid with the plane and how the colors vary along a strand of the saddle point braid as we traverse the interval $[2\pi k/\ell,2\pi (k+1)/\ell]$ between such intersection points.

\begin{figure}[h]
\centering
\includegraphics[height=5cm]{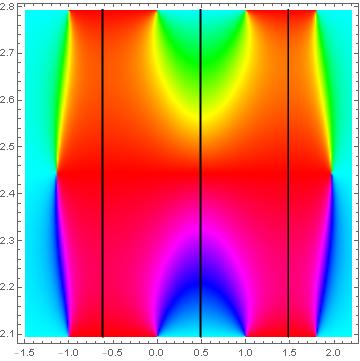}\qquad
\labellist
\Large
\pinlabel a) at -380 330
\pinlabel b) at -20 330
\pinlabel c) at -380 -40
\pinlabel d) at 10 -40
\tiny
\pinlabel $\text{Re}(u)$ at -50 -10
\pinlabel $t$ at -360 300
\pinlabel $\text{Re}(u)$ at 370 10
\pinlabel $t$ at 5 370
\pinlabel $\text{Re}(u)$ at -10 -375
\pinlabel $t$ at -365 -30
\pinlabel $\text{Re}(u)$ at 370 -375
\pinlabel $t$ at 25 -30
\endlabellist
\includegraphics[height=6cm]{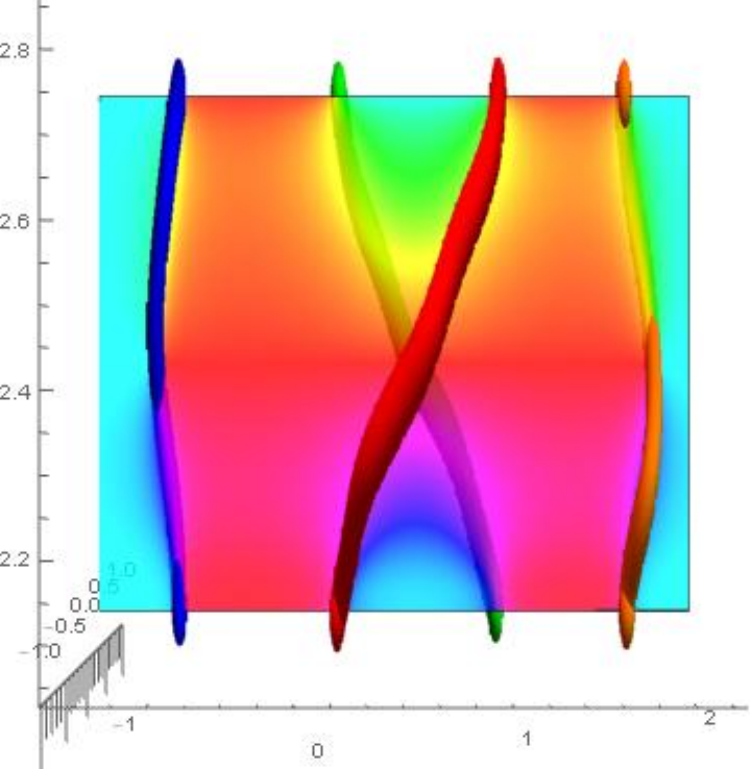}\\
\ \\
\includegraphics[height=5.5cm]{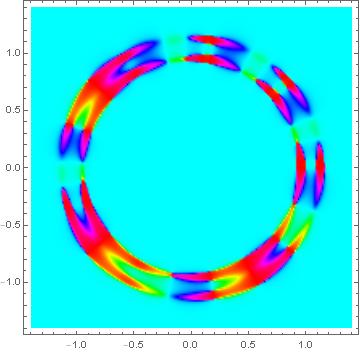}\qquad
\includegraphics[height=5.5cm]{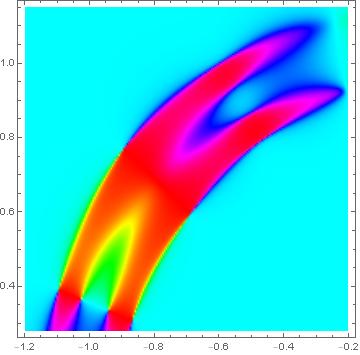}
\caption{a) The values of $\arg(g)$ on the $(\text{Im}(u)=0)$-plane in the interval $t\in[2\pi/3,8\pi/9]$. b) The roots of the $g$ in the interval $t\in[2\pi/3,8\pi/9]$ together with the values of $\arg(g)$ on the $(\text{Im}(u))$-plane. c) The values of $\arg(F)$ in the $(z=0)$-plane. d) The values of $\arg(F)$ in a region of the $(z=0)$-plane that contains the azimuthal interval $\varphi\in[2\pi/3,8\pi/9]$. \label{fig:72imu0_zoom}}
\end{figure}

Figure~\ref{fig:72imu0_zoom}a) shows a plot of $\arg(g_t)$ on the interval $[2\pi k/9,2\pi (k+1)/9]$ with $k=3$. Note that along the left-most critical point $(c_1,t)$ and along the right-most critical point $(c_3,t)$ the color is almost constant in this interval. It always stays in some red range. This means that the winding numbers of $\arg(g)(c_j,t):[2\pi k/9,2\pi (k+1)/9]\to S^1$ with $j=1$ and $j=3$ are both zero. Meanwhile, the colors on $(c_2,t)$ traverse the whole color wheel exactly once in a clockwise, i.e, "positive", direction, from cyan to blue, purple, red, yellow, green and finally cyan again. Thus the corresponding winding number is 1. In Figure~\ref{fig:72imu0_zoom}b) we see how this set of winding numbers generates a positive crossing $\sigma_2$ in this interval. 

It appears that there is a range of values for the parameters $a$ and $\mu$ where the solution to the paraxial wave equation with fixed behaviour on the $(z=0)$-plane behaves in a similar fashion as a complex-valued function that is holomorphic in $R+\rmi z$ with the same specified behaviour on the $(z=0)$-plane. In particular, values of the field on the annuli determine how the different roots twist around each other in each $2\pi/\ell$-interval of the azimuthal coordinate.

\begin{figure}[h]
\centering
\includegraphics[height=6cm]{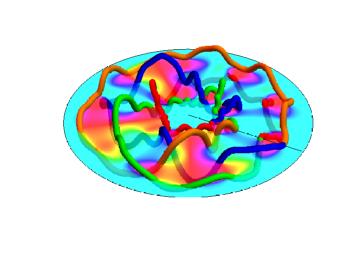}
\includegraphics[height=4cm]{optical72}
\caption{The zeros of $\Psi$ and $\arg(\Psi)$ on the $(z=0)$-plane. The zeros form the knot $7_2$.}
\end{figure}

Obviously not every complex-valued map in $\mathbb{R}^3$ behaves in this way. After all, the nodal topology of the polynomial field is different from the Gaussian beam with small width even if they are built from the same loop of complex polynomials $g_t$. Their zeros are mirror images (up to isotopy). We only claim that there is a choice of $a$ and $\mu$ such that the constructed Gaussian beam behaves like a holomorphic map. For other values, it does not and thus can have a different nodal topology.

We thus have two arguments why the two constructions produce the desired knot. The polynomial beam is sufficiently close to the constructed polynomial map $F$, so that it shares its nodal topology. The Gaussian beam behaves for some parameters like a holomorphic map and thus produces a knot from the $(z=0)$-slice in the same way that the zeros of $g_t$ are encoded in the $(\text{Im}(u)=0)$-slice. But why then do the two constructions not result in the same knot? Why do they produce two knots that are mirror images of each other instead?

The key to this question lies in the different substitutions and projection maps that we use in the constructions. The zeros of $g_t$ that form the input braid $B$ are in the space $\mathbb{C}\times[0,2\pi]\cong\mathbb{R}^2\times[0,2\pi]$, where $(\text{Re}(u),\text{Im}(u),t)$ forms a right-handed coordinate system. For every crossing we defined the overpassing strand to be the strand whose $\text{Im}(u)$-coordinate is smaller. This makes sense in so far as it is consistent with plots of the curves in this right-handed coordinate system, meaning that the overpassing strand lies above the undercrossing strand, compare for example Figure~\ref{fig:72imu0}c) and Figure~\ref{fig:braids}b).

In order to define the polynomial map $F$ we replace every instance of $\rme^{\rmi t}$ in $g_t$ by a complex variable $v$ and every $\rme^{-\rmi t}$ by its conjugate $\bar{v}$, and compose the resulting semiholomorphic polynomial with a stereographic projection. This means that $t$ is now identified with the azimuthal coordinate $\varphi$ in cylindrical coordinates in $\mathbb{R}^3$. The $\text{Re}(u)$-coordinate corresponds to the radial coordinate $R$ in cylindrical coordinates, in the sense that if $(x_1,x_2,x_3,x_4)$ and $(y_1,x_2,x_3,x_4)$ are two points on $S^3\subset \mathbb{C}^2\cong\mathbb{R}^4$ that only differ in their $\text{Re}(u)$-coordinate $x_1\neq y_1$, then their images $(R_1,\varphi_1,z_1)$ and $(R_2,\varphi_2,z_2)$ in cylindrical coordinates in $\mathbb{R}^3$ have the same azimuthal coordinate $\varphi_1=\varphi_2$ and $R_1<R_2$ if and only if $x_1<y_1$.

Similarly, the $\text{Im}(u)$-coordinate corresponds to the $z$-coordinate in cylindrical coordinates in $\mathbb{R}^3$. Note however that the coordinates in the right-handed coordinate system $(R,\varphi,z)$ have a different order, i.e., induce a different orientation, than their corresponding coordinates $(\text{Re}(u),\text{Im}(u),t)$. In other words, the map that sense a point $(a,b,c)\in\mathbb{R}^2\times[0,2\pi]$ in the $(\mathbb{Re}(u),\text{Im}(u),t)$-coordinates to $(a,c,b)\in\mathbb{R}_+\times[0,2\pi]\times\mathbb{R}$ in the $(R,\varphi,z)$-coordinate system is orientation-reversing. Before, we called at every crossing the strand with the smaller $\text{Im}(u)$-coordinate the overpassing strand. In $\mathbb{R}^3$, where the $\text{Im}(u)$-coordinate corresponds to the $z$-coordinate, this means that the strand with smaller $z$-coordinate is the overpassing strand. This is clearly the exact opposite of what we would expect for a reasonable definition of ``over''. The overpassing strand should lie above the other strand, i.e., it should have larger $z$-coordinate. This means that all crossings of the zeros of $F$ have exactly the opposite sign as the corresponding crossing in the desired knot, i.e., the optical vortex knot is the mirror image of the closure of the braid $B$. Since the polynomial beam inherits its nodal topology from $F$, the vortex knot is also the mirror image of the closure of $B$. 

The argument for the nodal topology of the Gaussian beam is independent of the nodal topology of $F$. As we traverse an interval $[2\pi k/\ell,2\pi(k+1)/\ell]$ in the azimuthal coordinate, roots that are close to each other twist around each other according to the winding number of $\arg(\Psi)$ along a path in the $(z=0)$-plane that lies between the roots. A positive winding number corresponds to positive total angular velocity, i.e., the motion of the roots in the $(R,z)$-plane is counterclockwise. This is the same as for the roots of the loop of complex polynomials $g_t$. Therefore, the crossings signs of the optical vortex knot are the same as for the braid closure, which means that it is the desired knot and the mirror image of the zeros of $F$ or the optical vortex knot of the polynomial beam.

We would like to emphasize that the relevant information about the knot that is stored in the $(z=0)$-plane is stable under perturbation, since it refers to topological quantities (winding numbers) associated to values of a smooth function on loops on the $(z=0)$-plane. That is, if the coefficients of $F(x,y,0)$ are approximated sufficiently well, the propagation of the approximated polynomial (times the Gaussian factor) also gives the desired optical vortex knot for approriate beam widths $w$. This is important, since in experiments it is not possible to prescribe the values of the field on the hologram plane $z=0$ exactly.

\section{The vortex knot table}\label{sec:table}
We apply our construction to all knots in the knot table up to 8 crossings. Since we are using braids and the minimal crossing number is not necessarily realised by a knot diagram in closed braid form, some braids have more than 8 crossings. The largest number of crossings that appears is 10. The largest number of strands is 5.

The loops $\gamma_j$, $j=1,2,\ldots,s$ and the base-point polynomial $p=g_0(u)$ can be chosen to only depend on $s$, the number of strands of the desired braid. In the construction of polynomial beams we used the following loops.
\begin{center}
\begin{tabular}{|c|c|c|}
\hline&&\\[-1em]
$s$ & $\gamma_j$ & $p=g_0$\\
\hline&&\\[-1em]
2 & $\gamma_2(\chi)=0.38 (\rme^{\rmi \chi}-1)$, $\chi\in[0,2\pi]$ & $g_0(u)=u(u+1)(u-1)$\\
\hline&&\\[-1em]
3 & $\gamma_2(\chi)=-1.5\rme^{\rmi \chi}+1.5$, $\chi\in[0,2\pi]$ & $g_0(u)=u(u-1.4)(u+1.4)$\\
 & $\gamma_3(\chi)=\rme^{\rmi \chi}-1$, $\chi\in[0,2\pi]$ & $\times(u-2.2)$ \\
\hline&&\\[-1em]
4 & $\gamma_2(\chi)=1.3(\rme^{\rmi \chi}-1)$, $\chi\in[0,2\pi]$ & $g_0(u)=u(u-1)(u+1)$  \\
& $\gamma_3(\chi)=-1.3(\rme^{\rmi \chi}-1)$, $\chi\in[0,2\pi]$ & $\times(u-2)(u+2)$\\
 & $\gamma_4(\chi)=2.2(\cos(\chi)-1)-1.5\rmi\sin(2\chi)$, $\chi\in[0,\pi]$ & \\
 & $2.2(\cos(\chi)-1)+1.5\rmi\sin(\chi)$, $\chi\in[\pi,2\pi]$, & \\
\hline&&\\[-1em]
5 & $\gamma_2(\chi)=-(\cos(\chi)-1)-\rmi\sin(\chi)$, $\chi\in[0,\pi]$ & $g_0(u)=u(u-0.7)(u+0.5)$\\
&$-(\cos(t)-1)+\rmi\sin(2\chi)$, $\chi\in[\pi,2\pi]$ & $\times(u-1.4)(u+1.5)(u+2.25)$ \\
 & $\gamma_3(\chi)=0.3(\cos(\chi)-1)+\rmi\sin(\chi)$, $\chi\in[0,2\pi]$ & \\
 & $\gamma_4(\chi)=-0.55(\cos(\chi)-1)-\rmi \sin(\chi)$, $\chi\in[0,2\pi]$ & \\
 & $\gamma_5(\chi)=1.4(\cos(\chi)-1)-1.4\rmi\sin(2\chi) $, $\chi\in[0,\pi]$ &\\
 &$1.4(\cos(\chi)-1)+1.4\rmi\sin(\chi)$, $\chi\in[\pi,2\pi]$ & \\
\hline
\end{tabular}
\end{center}

Here $s$ refers to the number of strands of the desired braid $B$. Recall that the degree of the complex polynomials $g_t$ is $s+1$, since we apply the algorithm by Bode-Hirasawa to the braid $\tilde{B}$, which has one more strand. Since this additional vertical strand does not interact with the rest of the braid, the generator $\sigma_1$ does not appear with either sign in the braid word of $\tilde{B}$ and thus the basic loop $\gamma_1$ is not needed for the construction of $\Gamma$.

The piecewise parametrisations used here for higher numbers of strands are easier to find and usually lead to a clearer separation between the strands of the braid. However, it takes longer to compute their Fourier transforms, which is one reason, beside the higher number of strands, why the construction of polynomial beams is computationally slightly more expensive than the construction of narrow Gaussian beams.

The table below lists the loops $\gamma_j$, $j=1,2,\ldots,s-1$ that were used in the construction of narrow Gaussian beams.
\begin{center}
\begin{tabular}{|c|c|c|}
\hline&&\\[-1em]
$s$ & $\gamma_j$ & $p=g_0$\\
\hline&&\\[-1em]
2 & $\gamma_1(\chi)=\tfrac{1}{4}(\rme^{\rmi \chi}-1)$, $\chi\in[0,2\pi]$ & $g_0(u)=u(u+1)$\\
\hline&&\\[-1em]
3 & $\gamma_1(\chi)=-\rme^{\rmi \chi}+1$, $\chi\in[0,2\pi]$ & $g_0(u)=u(u-1.4)(u+1.4)$\\
 & $\gamma_2(\chi)=\rme^{\rmi \chi}-1$, $\chi\in[0,2\pi]$ & \\
\hline&&\\[-1em]
4 & $\gamma_1(\chi)=0.53(\cos(\chi)-1)$ & $g_0(u)=u(u-1)(u+1)$  \\
& $+\frac{\rmi}{2}\left(0.55\sin(\chi)-\sin(\chi)^2\right.$ & $\times(u-1.8)$\\
 & $-0.2\cos(\chi)+0.2)$, $\chi\in[0,2\pi]$, & \\
 & $\gamma_2(\chi)=0.5(1-\rme^{\rmi  \chi})$, $\chi\in[0,2\pi]$, & \\
 & $\gamma_3(\chi)=0.375(\rme^{\rmi \chi}-1)$, $\chi\in[0,2\pi]$ & \\
\hline&&\\[-1em]
5 & $\gamma_1(\chi)=-2(\cos(\chi)-1)$ & $g_0(u)=u(u-1)(u+1)$\\
&$+\rmi (-\sin(\chi)^2-0.55\sin(\chi)$ & $\times(u-2)(u+2)$ \\
 & $-0.2\cos(\chi)+0.2)$, $\chi\in[0,2\pi]$ &  \\
 & $\gamma_2(\chi)=1.5(\cos(\chi)-1)+1.5\rmi\sin(\chi)$, $\chi\in[0,2\pi]$ & \\
 & $\gamma_3(\chi)=-1.5(\cos(\chi)-1)-1.5\rmi \sin(\chi)$, $\chi\in[0,2\pi]$ & \\
 & $\gamma_4(\chi)=2(\cos(\chi)-1) $ &\\
 &$+\rmi(-\sin(\chi)^2+0.55\sin(\chi)$ & \\
 & $-0.2\cos(\chi)+0.2)$, $\chi\in[0,2\pi]$ & \\
\hline
\end{tabular}
\end{center}



In the following table we list the knots together with information on the braid that we used and the values of the parameters that were used in order to produce the optical vortex knot. The parameter $m$ is the order of the Fourier approximation of the loop $\Gamma$. For example, for $m=12$ we used the Fourier coefficients of $\Gamma$ corresponding to $\rme^{-12\rmi t}, \rme^{-11\rmi t},\ldots,\rme^{11\rmi t},\rme^{12\rmi t}$.

$a$ is the scaling parameter described in the earlier sections, which is used in the construction of polynomial maps $F$, and $\mu$ is the parameter connected to the beam width, which only appears in the construction of Gaussian beams. The parameters $m$ and $a$ appear both in the construction of polynomial beams and in the construction of Gaussian beams. The index 1, e.g., in $m_1$, refers to the value of the parameter that was used in the construction of the Gaussian beam. The index 2, e.g., in $m_2$, refers to the value of the parameter that was used for the polynomial beam. 

We also include a column that specifies whether a given knot is a torus knot or a lemniscate knot. If the corresponding field is left empty, it means that it is neither. As we can see, most knots in the table are neither, which illustrates how far our constructions go beyond the known examples.

\begin{center}
\begin{tabular}{ |c|c|c|c|c|c|c|c|c| } 
 \hline&&&&&&&&\\[-1em]
 Knot & $s$ & braid word & $m_1$ & $a_1^{-1}$ & $\mu^{-1}$ & $m_2$ & $a_2^{-1}$ & Torus/Lemni.?\\
 \hline&&&&&& & &\\[-1em]
 $3_1$ & 2 & $\sigma_1^3$ & 3 & $6$ & 8 & 3 & 4 & Torus\\
 \hline&&&&&& & &\\[-1em]
 $4_1$ & 3 & $(\sigma_1\sigma_2^{-1})^2$ & 12 & $2$ & 1 & 12 & 4 & Lemniscate\\
 \hline&&&&&& & &\\[-1em]
 $5_1$ & 2 & $\sigma_1^5$ & 5 & $6$ & 8 & 5 & 4 & Torus\\
 \hline&&&&&& & &\\[-1em]
 $5_2$ & 3 & $\sigma_1^3\sigma_2\sigma_1^{-1}\sigma_2$ & 12 & 12 & 8 & 12 & 6 & \\
 \hline&&&&&& & &\\[-1em]
$6_1$ & 4 & $\sigma_1^2\sigma_2\sigma_1^{-1}\sigma_3^{-1}\sigma_2\sigma_3^{-1}$ & 20 & 6 & 8 & 20 & 8 & \\
 \hline&&&&&& & &\\[-1em]
$6_2$ & 3 & $\sigma_1^3\sigma_2^{-1}\sigma_1\sigma_2^{-1}$ & 6 & 6 & 8 & 6 & 6 & \\
 \hline&&&&&& & &\\[-1em]
$6_3$ & 3 & $\sigma_1^2\sigma_2^{-1}\sigma_1\sigma_2^{-1}$ & 12 & 6 & 8 & 12 & 6 & Lemniscate \\
 \hline&&&&&& & &\\[-1em]
$7_1$ & 2 & $\sigma_1^7$ & 7 & 6 & 8 & 7 & 4 & Torus \\
 \hline&&&&&& & &\\[-1em]
$7_2$ & 4 & $\sigma_1^3\sigma_2\sigma_1^{-1}\sigma_2\sigma_3\sigma_2^{-1}\sigma_3$ & 20 & 12 & 12 & 20 & 8 & \\
 \hline&&&&&& & &\\[-1em]
$7_3$ & 3 & $\sigma_1^5\sigma_2\sigma_1^{-1}\sigma_2$ & 12 & 6 & 8 & 20 & 8 & \\
 \hline&&&&&& & &\\[-1em]
$7_4$ & 4 & $\sigma_1^2\sigma_2\sigma_1^{-1}\sigma_2^2\sigma_3\sigma_2^{-1}\sigma_3$ & 20 & 8 & 12 & 20 & 8 & \\
 \hline&&&&&& & &\\[-1em]
$7_5$ & 3 & $\sigma_1^4\sigma_2\sigma_1^{-1}\sigma_2^2$ & 12 & 6 & 8 & 12 & 12 & \\
 \hline&&&&&& & &\\[-1em]
$7_6$ & 4 & $\sigma_1^2\sigma_2^{-1}\sigma_1\sigma_3\sigma_2^{-1}\sigma_3$ & 20 & 8 & 12 & 20 & 12 & \\
 \hline&&&&&& & &\\[-1em]
$7_7$ & 4 & $\sigma_1\sigma_2^{-1}\sigma_1\sigma_2^{-1}\sigma_3\sigma_2^{-1}\sigma_3$ & 20 & 12 & 8 & 20 & 12 & Lemniscate \\
 \hline&&&&&& & &\\[-1em]
$8_1$ & 5 & $\sigma_1^2\sigma_2\sigma_1^{-1}\sigma_2\sigma_3\sigma_2^{-1}\sigma_4^{-1}\sigma_3\sigma_4^{-1}$ & 25 & 40 & 16 & 40 & 8 & \\
 \hline&&&&& & &&\\[-1em]
$8_2$ & 3 & $\sigma_1^5\sigma_2^{-1}\sigma_1\sigma_2^{-1}$ & 12 & 6 & 8 & 12 & 6 & \\
 \hline&&&&& & &&\\[-1em]
$8_3$ & 5 & $\sigma_1^2\sigma_2\sigma_1^{-1}\sigma_3^{-1}\sigma_2\sigma_3^{-1}\sigma_4^{-1}\sigma_3\sigma_4^{-1}$ & 25 & 40 & 25 & 40 & 8 & \\
 \hline&&&&& & &&\\[-1em]
$8_4$ & 4 & $\sigma_1^3\sigma_2^{-1}\sigma_1\sigma_2^{-1}\sigma_3^{-1}\sigma_2\sigma_3^{-1}$ & 20 & 20 & 16 & 20 & 12 & \\
 \hline&&&&& & &&\\[-1em]
$8_5$ & 3 & $(\sigma_1^3\sigma_2^{-1})^2$ & 12 & 20 & 8 & 12 & 6 & \\
 \hline&&&&& & &&\\[-1em]
$8_6$ & 4 & $\sigma_1^4\sigma_2\sigma_1^{-1}\sigma_3^{-1}\sigma_2\sigma_3^{-1}$ & 20 & 30 & 16 & 20 & 12 & \\
 \hline&&&&& & &&\\[-1em]
$8_7$ & 3 & $\sigma_1^4\sigma_2^{-1}\sigma_1\sigma_2^{-2}$ & 12 & 20 & 8 & 12 & 6 & \\
 \hline&&&&& & &&\\[-1em]
$8_8$ & 4 & $\sigma_1^3\sigma_2\sigma_1^{-1}\sigma_3^{-1}\sigma_2\sigma_3^{-2}$ & 20 & 30 & 20 & 20 & 12 & \\
 \hline&&&&& & &&\\[-1em]
$8_9$ & 3 & $\sigma_1^3\sigma_2^{-1}\sigma_1\sigma_2^{-3}$ & 12 & 20 & 8 & 12 & 6 & Lemniscate \\
 \hline&&&&& & &&\\[-1em]
$8_{10}$ & 3 & $\sigma_1^3\sigma_2^{-1}\sigma_1^2\sigma_2^{-2}$ & 12 & 20 & 8 & 12 & 6 & \\
 \hline&&&&& & &&\\[-1em]
$8_{11}$ & 4 & $\sigma_1^2\sigma_2\sigma_1^{-1}\sigma_2^2\sigma_3^{-1}\sigma_2\sigma_3^{-1}$ & 20 & 30 & 40 & 20 & 12 & \\
 \hline&&&&& & &&\\[-1em]
$8_{12}$ & 5 & $(\sigma_1\sigma_2^{-1}\sigma_3\sigma_4^{-1})^2$ & 25 & 40 & 40 & 40 & 8 & Lemniscate \\
 \hline&&&&& & &&\\[-1em]
$8_{13}$ & 4 & $\sigma_1^2\sigma_2^{-1}\sigma_1\sigma_2^{-2}\sigma_3^{-1}\sigma_2\sigma_3^{-1}$ & 20 & 40 & 40 & 20 & 12 & \\
 \hline&&&&& & &&\\[-1em]
$8_{14}$ & 4 & $\sigma_1^3\sigma_2\sigma_1^{-1}\sigma_2\sigma_3^{-1}\sigma_2\sigma_3^{-1}$ & 20 & 40 & 40 & 20 & 12 & \\
 \hline&&&&& & &&\\[-1em]
$8_{15}$ & 4 & $\sigma_1^2\sigma_2^{-1}\sigma_1\sigma_3\sigma_2^3\sigma_3$ & 20 & 40 & 40 & 20 & 12 & \\
 \hline&&&&& & &&\\[-1em]
$8_{16}$ & 3 & $\sigma_1^2\sigma_2^{-1}\sigma_1^2\sigma_2^{-1}\sigma_1\sigma_2^{-1}$ & 12 & 20 & 20 & 12 & 6 & \\
 \hline&&&&& & &&\\[-1em]
$8_{17}$ & 3 & $\sigma_1^2\sigma_2^{-1}\sigma_1\sigma_2^{-1}\sigma_1\sigma_2^{-2}$ & 12 & 20 & 20 & 12 & 6 & \\
 \hline&&&&& & &&\\[-1em]
$8_{18}$ & 3 & $(\sigma_1\sigma_2^{-1})^4$ & 12 & 20 & 20 & 12 & 6 & Lemniscate \\
 \hline&&&&& & &&\\[-1em]
$8_{19}$ & 3 & $\sigma_1^3\sigma_2\sigma_1^{3}\sigma_2$ & 12 & 20 & 20 & 12 & 6 & Torus \\
 \hline&&&&& & &&\\[-1em]
$8_{20}$ & 3 & $\sigma_1^3\sigma_2^{-1}\sigma_1^{-3}\sigma_2^{-1}$ & 12 & 20 & 20 & 12 & 6 & \\
 \hline&&&&& & &&\\[-1em]
$8_{21}$ & 3 & $\sigma_1^3\sigma_2\sigma_1^{-2}\sigma_2^2$ & 12 & 20 & 20 & 12 & 6 & \\
 \hline
\end{tabular}
\end{center}

The values are not necessarily optimal. It is for example conceivable that the knot $7_2$ can be reproduced with larger values of $a$.

\begin{figure}
\centering
\includegraphics[height=3cm]{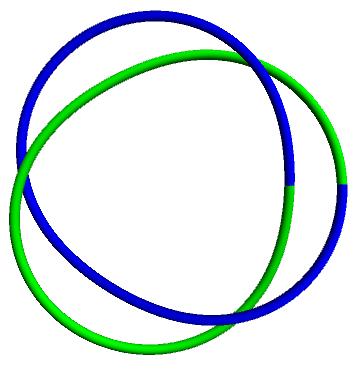}
\includegraphics[height=3cm]{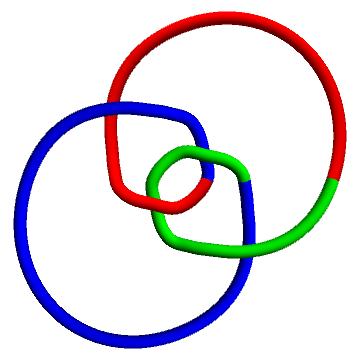}
\includegraphics[height=3cm]{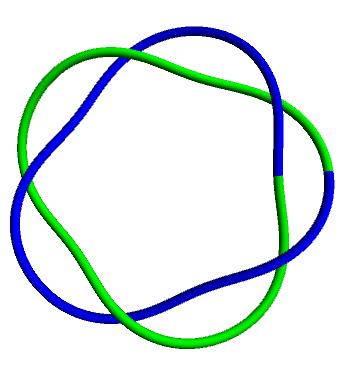}
\includegraphics[height=3cm]{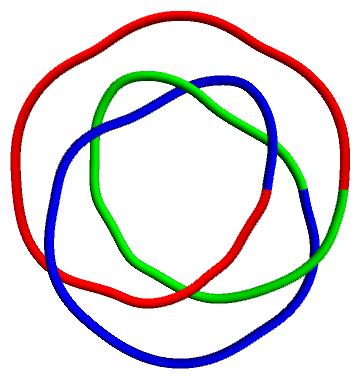}
\includegraphics[height=3cm]{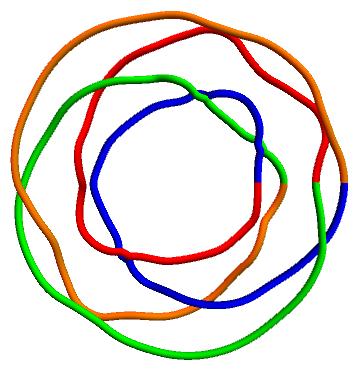}
\includegraphics[height=3cm]{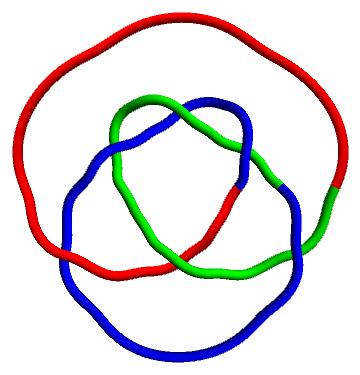}
\includegraphics[height=3cm]{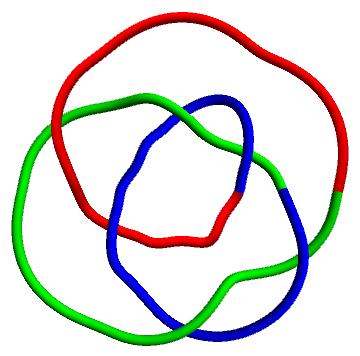}
\includegraphics[height=3cm]{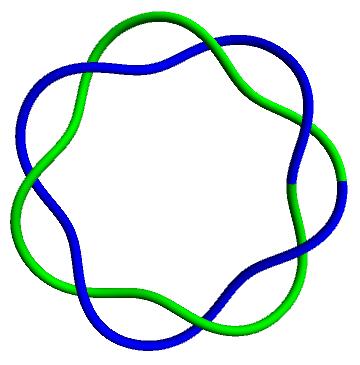}
\includegraphics[height=3cm]{optical72}
\includegraphics[height=3cm]{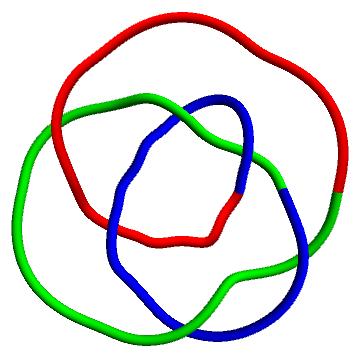}
\includegraphics[height=3cm]{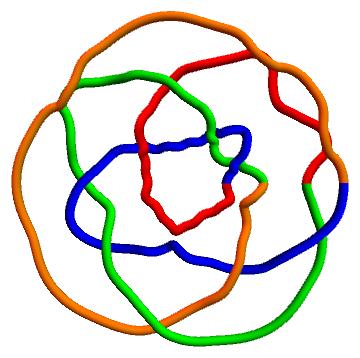}
\includegraphics[height=3cm]{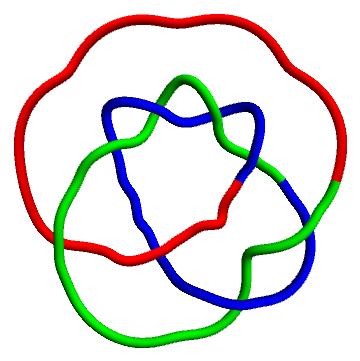}
\includegraphics[height=3cm]{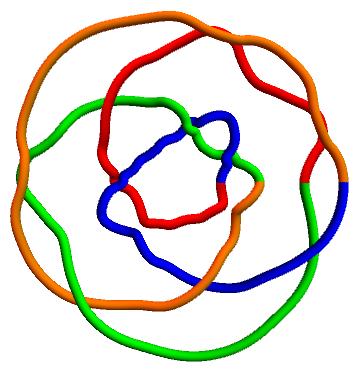}
\includegraphics[height=3cm]{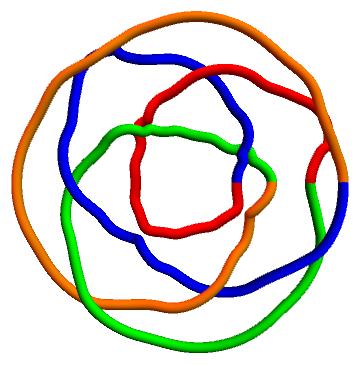}
\caption{The table of optical vortex knots of up to 7 crossings. The curves are the traced zeros of the constructed optical fields for each knot type.\label{fig:knottable}}
\end{figure}

\begin{figure}
\centering
\includegraphics[height=3cm]{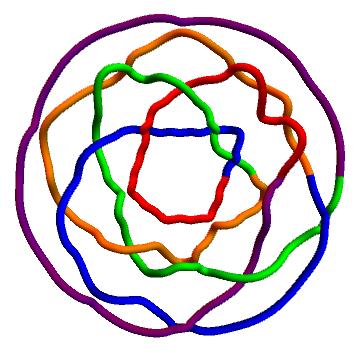}
\includegraphics[height=3cm]{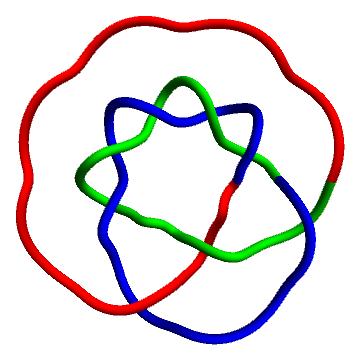}
\includegraphics[height=3cm]{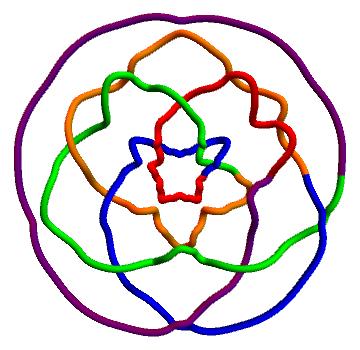}
\includegraphics[height=3cm]{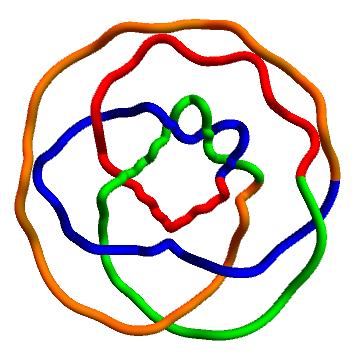}
\includegraphics[height=3cm]{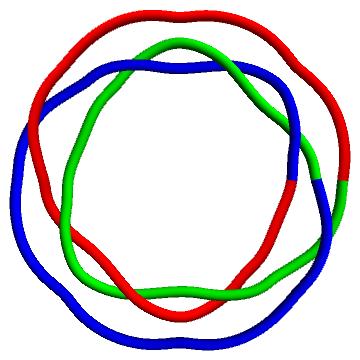}
\includegraphics[height=3cm]{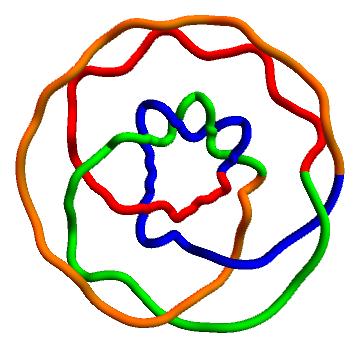}
\includegraphics[height=3cm]{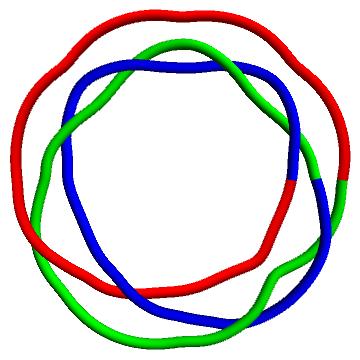}
\includegraphics[height=3cm]{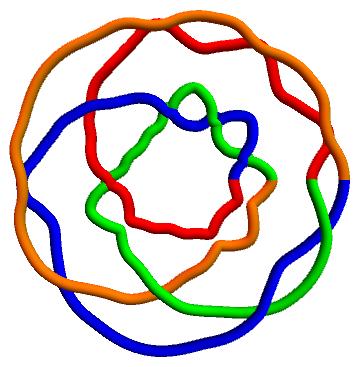}
\includegraphics[height=3cm]{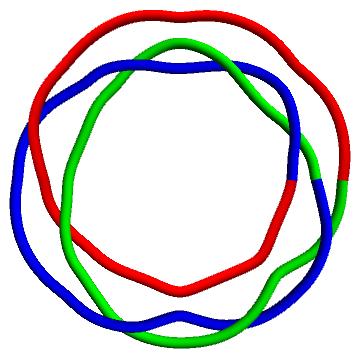}
\includegraphics[height=3cm]{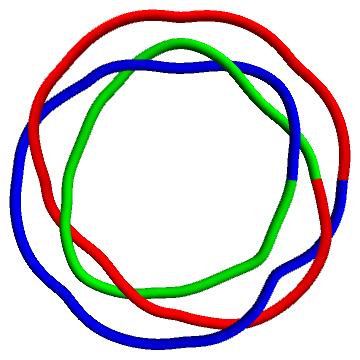}
\includegraphics[height=3cm]{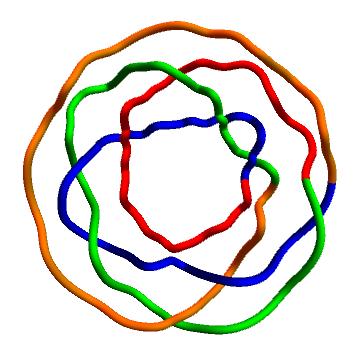}
\includegraphics[height=3cm]{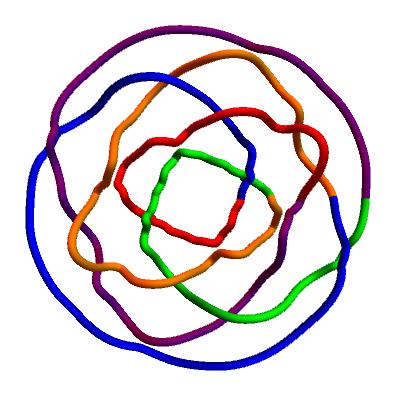}
\includegraphics[height=3cm]{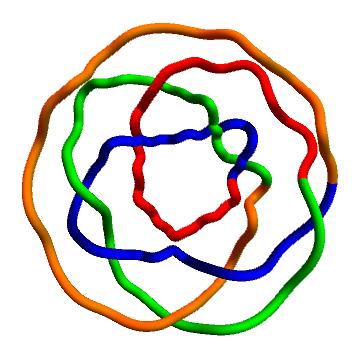}
\includegraphics[height=3cm]{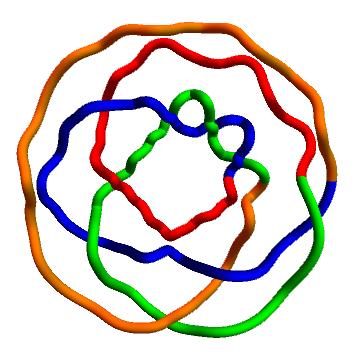}
\includegraphics[height=3cm]{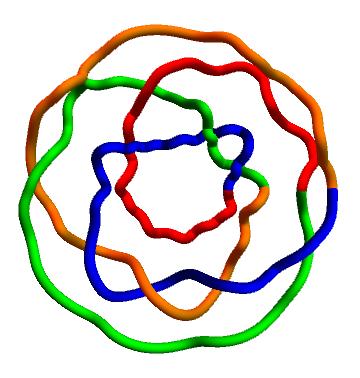}
\includegraphics[height=3cm]{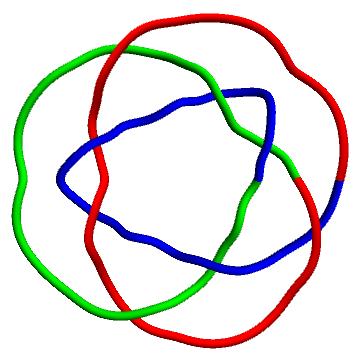}
\includegraphics[height=3cm]{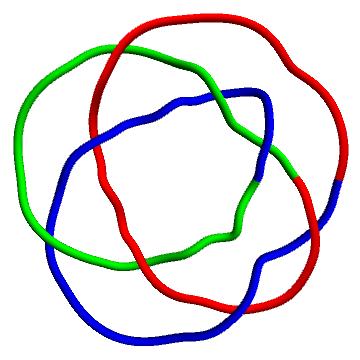}
\includegraphics[height=3cm]{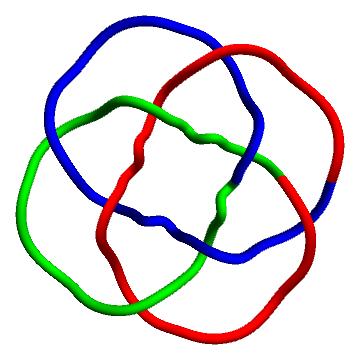}
\includegraphics[height=3cm]{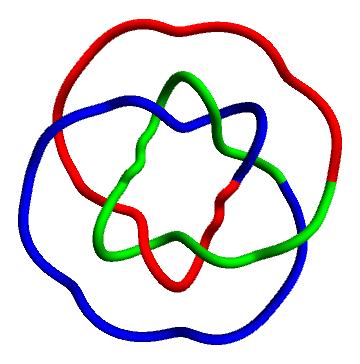}
\includegraphics[height=3cm]{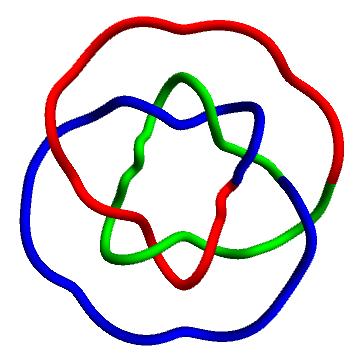}
\includegraphics[height=3cm]{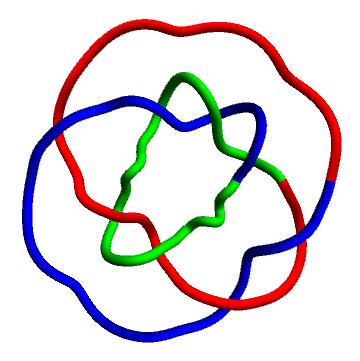}
\caption{The table of optical vortex knots of 8 crossings. The curves are the traced zeros of the constructed optical fields for each knot type.\label{fig:knottable2}}
\end{figure}

Below, in Figures~\ref{fig:knottable} and \ref{fig:knottable2}, we plot the zeros obtained from the construction of narrow Gaussian beams. As in the example of the knot $7_2$, the curves are found by tracking the zeros of the field numerically through 840 fields of constant azimuthal coordinate $\varphi=2\pi j/840$, $j=0,1,2,\ldots,839$. This means that in principle there could be other components of the nodal set and in general it looks like this is the case.

The different colours along the knot correspond to different strands of the braid. Very small values of the parameter $a$ imply that the knot lies in a small tubular neighbourhood of the planar circle $x^2+y^2=\mu^{-2}$, $z=0$. We display the knots as having roughly the same size by scaling the curves appropriately. We also scale them so that the strands are separated enough to distinguish which strand is passing over the other at a crossing. In this way, we can read off the braid word from the plots. These are not necessarily identical to the input braid word, but equivalent via some Reidemeister moves of the second type.


For each of these knots two mathematica files that contain the construction of the polynomial beam and of the Gaussian beam are published on the author's webpage \cite{webpage}. As we can see we have successfully constructed any knot type of up to 8 crossings as an optical vortex knot, once in a polynomial beam and once in a narrow Gaussian beam. The creation of these fields and the subsequent plot of the knot takes only few minutes for braids of 3. For one of the most complicated knots (5 strands plus one more vertical strand in the construction of polynomial beams) it takes several hours on a standard laptop computer. 

\section{Some comments on quantum vortex knots}\label{sec:quantum}

Note that the paraxial wave equation is a 2+1-dimensional Schrödinger equation without potential. Part of the original motivation for this project was to understand Dennis's propagation technique better, so that it may be applied to the construction of solutions to 3+1-dimensional Schrödinger equations, i.e., quantum wavefunctions whose zeros are surfaces that describe given motions of quantum vortex knots. The Gross-Pitaevski equation, which models Bose-Einstein condensates is an example of a non-linear Schrödinger equation. In this case, the general idea of propagating a slice of a polynomial map could still work. However, since the equation is non-linear, this propagation cannot be realized by simply propagating each monomial and summing the resulting fields.

In analogy with quantum vortex knots, our solutions of the paraxial wave equation can be interpreted as wave functions in 2 dimensions, whose vortices (points) trace out a given knot. The dynamics of point vortices in 2-dimensional fluids and superfluids are of interest in its own right \cite{aref}, but controlling the motion of quasiparticles in 2-dimensions so that they form knots has also become important in certain areas of topological quantum computing \cite{tqc}.

There are still many unanswered questions regarding the 3+1-dimensional setting. Given an initial wavefunction at $t=0$ whose zeros form a knot, it is extremely difficult to predict what will happen to the vortex knot. It is known theoretically that any cobordism between knots can be realised by the zeros of a quantum wavefunction (for example in the harmonic oscillator or in the Gross-Pitaevski equation \cite{daniel}). In other words, any type of motion and of topology change is possible. However, finding an initial configuration for a desired motion is an incredibly subtle endeavour and it took a long time until an example was discovered where a Hopf link changes into a more complicated knot: the trefoil knot \cite{simone}. 

As in the lower-dimensional case, it is well-known that every embedded surface in $\mathbb{R}^4$ can be realized (up to isotopy) as the zeros of a polynomial map $F:\mathbb{R}^4\to\mathbb{R}^2\cong\mathbb{C}$. There is a construction of complex-valued polynomials in four real variables, whose zeros form embedded (linked) tori \cite{kamada}, which describe loops moving around each other and through each other without undergoing any reconnection events. However, more complicated surfaces, that is surfaces with different genus (describing motions with reconnection events), seem to require approximation techniques instead of interpolation techniques, which were used for knotted tori. This means that the resulting polynomials can have very high degrees, while both for knots and tori there exist proven upper bounds on the topological degrees in terms of topological data of the knot or embedded torus in question \cite{bodepoly, kamada}. For linear Schrödinger type equations Dennis's propagation technique should directly generalize to 3+1 dimensions. Instead of restricting the polynomial to the ($z=0$)-plane and using this as a boundary condition, we take the polynomial in four variables and restrict it to the $x_4=0$-hyperplane. However, (even in the lower-dimensional case) it is not understood why the propagation should result in the same topological type of zeros. 

In order to find a polynomial map $F:\mathbb{R}^4\to\mathbb{R}^2$ whose zeros form a given surface $S$ in $\mathbb{R}^4$, we may consider a direct analogue of the construction by Bode-Hirasawa. There is a notion of a surface braid and as in the lower-dimensional case, any embedded, compact, oriented surface can be realized as the closure of such a surface braid \cite{kamadabraid, kamadabook}. We take $p$ once again to be a complex polynomial of degree $s$ with distinct real roots and distinct critical values. Let $\Gamma:S^2\to\mathbb{C}$ and $g_{\varphi,\theta}(u)=p(u)-\Gamma(\varphi,\theta)$, where $\varphi$ and $\theta$ are the spherical coordinates on $S^2$. In general, the zeros of $g$ form a braided surface in $\mathbb{C}\times S^2$. Choosing a scaling parameter $a$ sufficiently small means that the complex coordinate of the zeros is close to 0.

For knots we approximated the constant term (which was the only coefficient of $g_t$ that depended on $t$) by trigonometric polynomials and made the replacement $\rme^{\rmi t}\to v$, $\rme^{-\rmi t}\to\bar{v}$. Likewise, we can approximate $\Gamma(\varphi,\theta)$ by spherical harmonics and and replace each summand by the unique homogeneous polynomial in three variables that restricts to the monomial on $S^2$. In this way we obtain a polynomial in 3 real variables that agrees with $\Gamma$ on $S^2\subset\mathbb{R}^3$. We call the difference of $p$ and this polynomial in 3 real variables $f:\mathbb{C}\times\mathbb{R}^3\to\mathbb{C}$, which is a complex-valued polynomial in the complex variable $u$ and the three real variables. Writing $u$ as a pair of real variables, we have a polynomial map $f:\mathbb{R}^5\to\mathbb{R}^2$. The same arguments as in \cite{bodepoly} should imply that if $a$ is chosen sufficiently small, then the zeros of $f$ in $S^4\subset\mathbb{R}^5$ are the desired surface $S$, in the sense that a composition of $f$ with a stereographic projection map between $S^4$ and $\mathbb{R}^4\cup\{\infty\}$ produces a polynomial $F:\mathbb{R}^4\to\mathbb{R}^2$ whose zeros are the closure of the braided surface that is given by the zeros of $g_{\varphi,\theta}$.

In practice, this is a lot more complicated than the lower-dimensional setting and the propagation will be even more difficult. However, the theoretical arguments from our discussion on knots should still apply. By construction the constant term of $g$ is the only coefficient that depends on $\varphi$ and $\theta$. This means that the zeros of the derivative of $g$ with respect to $u$ do not depend on $\varphi$ and $\theta$ and are thus the critical points of $p$: the real numbers $c_j$, $j=1,2,\ldots,s-1$. In particular, the critical surfaces $(c_j,\varphi,\theta)$, $j=1,2,\ldots,s-1$, $\varphi\in[0,2\pi]$, $\theta\in[0,\pi]$, lie in the three-dimensional hyperplane $\text{Im}(u)=0$. Under stereographic projection this hyperplane gets mapped to the $(x_4=0)$-hyperplane. So the values that the polynomial $F$ takes on the critical surfaces (which should encode the topology of $S$ in some sense) can be found in the hyperplane that we use to propagate the field, leaving us hopeful that the resulting solution $\Psi$ of the 3+1-Schrödinger equation reproduces the surface of the polynomial map, that is, (up to isotopy and mirror images) the zeros of the field $\Psi$ form the desired surface $S$.

There are still various gaps to fill in these arguments and several details to work out, but at least conceptually, the construction of optical vortex knots provides a promising blueprint for the construction of embedded surfaces as zeros of 3+1-dimensional Schrödinger equations, which describe the time-evolution of quantum vortex knots.


\begin{thebibliography}{99}
\bibitem{alexander} J. W. Alexander. \textit{A lemma on systems of knotted curves}. Proceedings of the National Academy of Sciences of the United States of America \textbf{9}, 3 (1923), 93--95.
\bibitem{aref} H. Aref. \textit{Point vortex dynamics: A classical mathematics playground}. Journal of Mathematical Physics \textbf{28} (2007), 065401.
\bibitem{knotatlas} D. Bar-Natan \textit{The Knot Atlas}. \textit{http://katlas.org/}. (Accessed June 16th, 2024.) 
\bibitem{berry} M. V. Berry. \textit{Knotted zeros in the quantum states of hydrogen}. Found. Phys. \textbf{31} (2001), 659--667.
\bibitem{bd1} M. V. Berry and M. R. Dennis. \textit{Knotted and linked phase singularities in monochromatic waves}. Proc. R. Soc. Lond. A \textbf{457} (2001), 2251--2263.
\bibitem{bd2} M. V. Berry and M. R. Dennis. \textit{Knotting and unknotting of phase singularities: Helmholtz waves, paraxial waves and waves in 2+1 spacetime}. J. Phys. A \textbf{34} (2001), 8877.
\bibitem{RAG} Bochnak J, Coste M, Roy M-F. \textit{Real algebraic geometry} (Springer Verlag, Berlin, Heidelberg, 1998).
\bibitem{phd} B. Bode. \textit{Knotted fields and real algebraic links}. PhD thesis. University of Bristol (2018).
\bibitem{quasi} B. Bode. \textit{Quasipositive links and electromagnetism}. Topology and its Applications \textbf{301} (2021), 107495.
\bibitem{bodeem} B. Bode. \textit{Stable knots and links in electromagnetic fields}. Comm. Math. Phys. \textbf{387} (2021), 1757–1770.
\bibitem{webpage} B. Bode. \textit{Personal webpage}. \textit{https://sites.google.com/view/benjaminbode/research}. (Accessed June 18th, 2024.)
\bibitem{bodepoly} B. Bode and M. R. Dennis. \textit{Constructing a polynomial whose nodal set is any prescribed knot or link}. Journal of Knot Theory and its Ramifications \textbf{28} (2019), 1850082.
\bibitem{lemniscate} B. Bode, M. R. Dennis, D. Foster and R. P. King. \textit{Knotted fields and explicit fibrations for lemniscate knots}. Proc. R. Soc. A \textbf{473} (2017), 20160829.
\bibitem{mikami} B. Bode and M. Hirasawa. \textit{Saddle point braids of braided fibrations and pseudo-fibrations}. Research in the Mathematical Sciences (to appear).
\bibitem{kamada} B. Bode and S. Kamada. \textit{Knotted surfaces as vanishing sets of polynomials}. Journal of the Mathematical Society of Japan \textbf{73}, no. 4 (2021), 1289--1322.
\bibitem{holo} R. Bousso. \textit{The holographic principle}. Reviews of Modern Physics \textbf{74}, 3 (2002), 825--874.
\bibitem{brauner} K. Brauner. \textit{Zur Geometrie der Funktionen zweier komplexen Ver{\"a}nderlichen II, III, IV}. Abh. Math. Sem. Hamburg {\bf6} (1928) 8--54.
\bibitem{52} M. R. Dennis and B. Bode. \textit{Constructing a polynomial whose nodal set is the three-twist knot $5_2$}. Journal of Physics A \textbf{50} (2017), 265204.
\bibitem{marknature} M. R. Dennis, R. P. King, B. Jack, K. O'Holleran and M. J. Padgett. \textit{Isolated optical vortex knots}. Nat. Phys.\textbf{6} (2010). 118--121.
\bibitem{daniel} A. Enciso and D. Peralta-Salas. \textit{Approximation theorems for the Schrödinger equation and quantum vortex reconnection}. Comm. Math. Phys. \textbf{387} (2021), 1111-1149.
\bibitem{small} I. Herrera, A. Mojica-Casique and P. A. Quinto-Su. \textit{Experimental realization of a wavelength-sized optical-vortex knot}. Physical Review Applied \textbf{17} (2022), 064026.
\bibitem{kamadabraid} S. Kamada. \textit{A characterization of groups of closed orientable surfaces in 4-space}. Topology \textbf{33}, 1 (1994), 113--122.
\bibitem{kamadabook} S. Kamada. \textit{Braid and knot theory in dimension four}. American Mathematical Society, Providence, RI (2002).
\bibitem{kamien} R. D. Kamien and R. A. Mosna. \textit{The topology of dislocations in smectic liquid crystals}. New J. Phys. \textbf{18} (2016), 053012.
\bibitem{braids} C. Kassel and V. Turaev. \textit{Braid groups} (Springer-Verlag New York, 2008).
\bibitem{irvine} D. Kleckner and W. T. M. Irvine. \textit{Creation and dynamics of knotted vortices}. Nature Phys. \textbf{9} (2013), 253--258.
\bibitem{coding} L.-J. Kong, W. Zhang, P. Li, X. Guo, J. Zhang, F. Zhang, J. Zhao and X. Zhang. \textit{High capacity topological coding based on nested vortex knots and links}. Nature Communications \textbf{13} (2022), article no. 2705.
\bibitem{info} H. Larocque, A. D'Errico, M. F. Ferrer-Garcia, A. Carmi, E. Cohen and E. Karimi. \textit{Optical framed knots as information carriers}. Nature Communications \textbf{11} (2020), article no. 5119.
\bibitem{leach} J. Leach ,M. R. Dennis, J. Courtial and M. J. Padgett. \textit{Vortex knots in light}. New Journal of Physics \textbf{7} (205), 55.
\bibitem{cascades} X. Liu and R. Ricca. \textit{Knots cascade detected by a monotonically decreasing sequence of values}. Scientific Reports \textbf{6} (2016), 24118.
\bibitem{ma} Machon T, Alexander GP. 2014. Knotted defects in nematic liquid crystals. \textit{Phys. Rev. Lett.} \textbf{113}, 027801.
\bibitem{milnor} J. W. Milnor, \textit{Singular points of complex hypersurfaces} (Princeton University Press, 1968).
\bibitem{moffatt:1969degree} H.K Moffatt. \textit{The degree of knottedness of tangled vortex lines}. J. Fluid. Mech. \textbf{35} (1969), 117--129.
\bibitem{proment} D. Proment, M. Onorato and C. F. Barenghi, Vortex knots in a Bose-Einstein condensate, {\it Phys. Rev. E} {\bf 85}, 3 (2012) 036306.
\bibitem{bessel} T. Radozycki. \textit{Knotted nodal lines in superpositions of Bessel-Gaussian light beams}. Phys. Rev. A \textbf{103} (2021), 013509.
\bibitem{reidemeister} K. Reidemeister. \textit{Elementare Begründung der Knotentheorie}. Abhandlungen aus dem Mathematischen Seminar der Universität Hamburg \textbf{5} (1927), 24--32.
\bibitem{tqc} A. Stern and N. H. Lindner. \textit{Topological quantum computation-- From basic concepts to first experiments}. Science \textbf{339}, Issue 6124 (2013), 1179--1184. 
\bibitem{dna} R. Stolz, M. Yoshida, R. Brasher, M. Flanner, K. Ishihara, D. J. Sherratt, K. Shimokawa and M. Vazquez. \textit{Pathways of DNA unlinking: A story of stepwise simplification}. Scientific Reports \textbf{7} (2017), 12420.
\bibitem{dani} D. Sugic. \textit{Unravelling the dark focus of light: a study of knotted optical singularities}. PhD thesis. University of Bristol (2019).
\bibitem{susskind} L. Susskind. \textit{The world as a hologram}. Journal of Mathematical Physics \textbf{36}, 11 (1995), 6377--6396.
\bibitem{sutcliffe} Sutcliffe P. 2007. Knots in the Skyrme-Faddeev model. \textit{Proc. R. Soc. A} \textbf{463}, 3001--3020.
\bibitem{photon} S. J. Tempone-Wiltshire, S. P. Johnstone and K. Helmerson. \textit{Optical vortex knots -- one photon at a time}. Scientific Reports \textbf{6} (2016), 24463.
\bibitem{kelvin} W. Thomson (Lord Kelvin), \textit{On vortex atoms}. Proceedings of the Royal Society of Edinburgh \textbf{VI} (1867), 94--105.
\bibitem{colloid} U. Tkalec , M. Ravnik, S. Čopar, S. Žumer and Igor Muševič. \textit{Reconfigurable knots and links in chiral nematic colloids}. Science \textbf{333}, 6038 (2011), 62--65. 
\bibitem{ultrasmall} L. Wang, W. Zhang, H. Yin and X. Zhang. \textit{Ultrasmall optical vortex knots generated by spin-selective metasurface holograms}. Advanced Optical Materials \textbf{7}, 10 (2019), 1900263.
\bibitem{acoustic} H. Zhang, W. Zhang, Y. Liao, X. Zhou, J. Li, G Hu and X. Zhang. \textit{Creation of acoustic vortex knots}. Nature Communications \textbf{11} (2020), 3956.
\bibitem{charge} J. Zhong, S. Liu, X. Guo, P. Li, B. Wei, L. Han, S. Qi and J. Zhao. \textit{Observation of optical vortex knots and links associated with topological charge }. Optics Express \textbf{29}, Issue 23 (2021), 38849-38857.
\bibitem{simone} S. Zuccher and R. L. Ricca. \textit{Creation of quantum knots and links driven by minimal surfaces}. J. Fluid. Mech. \textbf{942} (2022), A 8.


\end{thebibliography}
\end{document}